\def\version{{\tiny Version \today (typeset: \today)}}
\def\version{}

\documentclass[reqno,twoside,11pt]{amsart}
\usepackage{cite}
\usepackage{amsmath}
\usepackage{amsfonts}
\usepackage{amssymb}
\usepackage{epsfig}
\usepackage{verbatim}
\usepackage{color}
\usepackage{bbm}
\IfFileExists{epsf.def}{\input epsf.def}{\usepackage{epsf}}

\IfFileExists{mathptmx.sty}{\usepackage{mathptmx}}{\usepackage{mathpazo}}
\usepackage{mathrsfs}

\DeclareFontFamily{OT1}{eusb}{} \DeclareFontShape{OT1}{eusb}{m}{n} {<5> <6> <7> <8> <9> <10> <11> <12> <14.4> eusb10}{}
\DeclareMathAlphabet{\eusb}{OT1}{eusb}{m}{n}

\DeclareFontFamily{OT1}{eusm}{} \DeclareFontShape{OT1}{eusm}{m}{n} {<5> <6> <7> <8> <9> <10> <11> <12> <14.4> eusm10}{}
\DeclareMathAlphabet{\eusm}{OT1}{eusm}{m}{n}

\DeclareFontFamily{OT1}{eufm}{} \DeclareFontShape{OT1}{eufm}{m}{n} {<5> <6> <7> <8> <9> <10> <11> <12> <14.4> eufm10}{}
\DeclareMathAlphabet{\mathfrak}{OT1}{eufm}{m}{n}

\DeclareFontFamily{OT1}{fraktura}{}
\DeclareFontShape{OT1}{fraktura}{m}{n} {<5> <6> <7> <8> <9> <10> <11> <12> <13> <14.4> [1.1] eufm10}{}
\DeclareMathAlphabet{\fraktura}{OT1}{fraktura}{m}{n}

\DeclareFontFamily{OT1}{cmfi}{} \DeclareFontShape{OT1}{cmfi}{m}{n} {<5> <6> <7> <8> <9> <10> <11> <12> <13> <14.4> [0.9] cmfi10}{}
\DeclareMathAlphabet{\cmfi}{OT1}{cmfi}{b}{n}

\DeclareFontFamily{OT1}{cmss}{} \DeclareFontShape{OT1}{cmss}{m}{n} {<5> <6> <7> <8> <9> <10> <11> <12> <13> <14.4> cmss10}{}
\DeclareMathAlphabet{\cmss}{OT1}{cmss}{m}{n}

\newtheoremstyle{thm}{1.5ex}{1.5ex}{\itshape\rmfamily}{} {\bfseries\rmfamily}{}{2ex}{}

\newtheoremstyle{def}{1.5ex}{1.5ex}{\rmfamily\sl}{} {\bfseries\rmfamily}{}{2ex}{}

\newtheoremstyle{rem}{1.3ex}{1.3ex}{\rmfamily}{} {\bfseries\rmfamily}{}{2ex}{}

\newtheoremstyle{ass}{1.5ex}{1.5ex}{\rmfamily\sl}{} {\bfseries\rmfamily}{}{2ex}{}

\newenvironment{proofsect}[1] {\vskip0.1cm\noindent{\rmfamily\itshape#1.}}{\qed\vspace{0.15cm}}

\theoremstyle{thm}
\newtheorem{theorem}{Theorem}[section]
\newtheorem{lemma}[theorem]{Lemma}
\newtheorem{proposition}[theorem]{Proposition}

\newtheorem*{Main Theorem}{Main Theorem.}

\newtheorem{conjecture}[theorem]{Conjecture}

\newtheoremstyle{named}{}{}{\itshape}{}{\bfseries}{}{.5em}{\thmnote{#3}}
\theoremstyle{named}

\theoremstyle{def}
\newtheorem{definition}[theorem]{Definition}

\theoremstyle{rem}
\newtheorem{remark}[theorem]{{Remark}}

\numberwithin{equation}{section}

\renewcommand{\section}{\secdef\sct\sect}
\newcommand{\sct}[2][default]{\refstepcounter{section}
\addcontentsline{toc}{section}
{{\tocsection {}{\thesection}{\!\!\!\!#1\dotfill}}{}}
\vspace{0.7cm}
\centerline{ 
\scshape\arabic{section}.\ #1} \nopagebreak \vspace{0.2cm}}
\newcommand{\sect}[1]{
\vspace{0.4cm} \centerline{\large\scshape\rmfamily #1}
\vspace{0.2cm}}

\renewcommand{\subsection}{\secdef\subsct\sbsect}
\newcommand{\subsct}[2][default]{\refstepcounter{subsection}
\addcontentsline{toc}{subsection}
{{\tocsection{\!\!}{\hspace{1.2em}\thesubsection}{\!\!\!\!#1\dotfill}}{}}
\nopagebreak\vspace{0.45\baselineskip} {\flushleft\bf
\arabic{section}.\arabic{subsection}~\bf #1.~}
\\*[3mm]\noindent
\nopagebreak}
\newcommand{\sbsect}[1]{
\vspace{0.1cm}\noindent
\textbf{#1.~}\vspace{0.1cm}}

\renewcommand{\subsubsection}{%
\secdef \subsubsect\sbsbsect}
\newcommand{\subsubsect}[2][default]{%
\refstepcounter{subsubsection} 
\addcontentsline{toc}{subsubsection}{{\tocsection{\!\!}
{\hspace{3.05em}\thesubsubsection}{\!\!\!\!#1\dotfill}}{}}
\nopagebreak
\vspace{0.15\baselineskip} \nopagebreak {\flushleft\rmfamily
\itshape\arabic{section}.\arabic{subsection}.\arabic{subsubsection}
\ \rmfamily #1\/.}\ }
\newcommand{\sbsbsect}[1]{\vspace{0.1cm}\noindent
\rmfamily \itshape
\arabic{section}.\arabic{subsection}.\arabic{subsubsection} \
\sffamily #1\/.\ }

\renewcommand{\caption}[1]{%
\vglue0.5cm
\refstepcounter{figure}
\begin{minipage}{0.9\textwidth}\small {\sc Figure~\thefigure. }#1\end{minipage}}

\newcommand{\dist}{\operatorname{dist}}
\newcommand{\supp}{\operatorname{supp}}
\newcommand{\diam}{\operatorname{diam}}

\newcommand{\Int}{{\text{\rm Int}}}

\newcommand{\textd}{\text{\rm d}\mkern0.5mu}

\newcommand{\texte}{\text{\rm  e}\mkern0.7mu}

\newcommand{\1}{{1\mkern-4.5mu\textrm{l}}}
\renewcommand{\1}{\ind}

\newcommand{\BB}{{\mathbb{B}}}

\newcommand{\BE}{{\mathbb{E}}}

\newcommand{\BH}{{\mathbb{H}}}

\newcommand{\BN}{{\mathbb{N}}}

\newcommand{\BP}{{\mathbb{P}}}

\newcommand{\BR}{{\mathbb{R}}}

\newcommand{\BZ}{{\mathbb{Z}}}

\newcommand{\CF}{{\mathcal{F}}}

\newcommand{\CM}{{\mathcal{M}}}

\newcommand{\CP}{{\mathcal{P}}}

\newcommand{\CU}{{\mathcal{U}}}

\newcommand{\ind}{{\mathbbm{1}}}

\newcommand{\bae}{\begin{equation}\begin{aligned}}
\newcommand{\eae}{\end{aligned}\end{equation}}

\newcommand{\Z}{\mathbb{Z}}

\DeclareFontFamily{OML}{rsfs}{\skewchar\font'177}
\DeclareFontShape{OML}{rsfs}{m}{n}{ <5> <6> rsfs5 <7> <8> <9>
rsfs7 <10> <10.95> <12> <14.4> <17.28> <20.74> <24.88> rsfs10 }{}
\DeclareMathAlphabet{\mathfs}{OML}{rsfs}{m}{n}

\newcommand{\twoeqref}[2]{(\ref{#1}--\ref{#2})}

\newcommand{\cc}{{\text{\rm c}}}

\def\myffrac#1#2 in #3{\raise 2.6pt\hbox{$#3 #1$}\mkern-1.5mu\raise 0.8pt\hbox{$#3/$}\mkern-1.1mu\lower 1.5pt\hbox{$#3 #2$}}
\newcommand{\ffrac}[2]{\mathchoice%
	{\myffrac{#1}{#2} in \scriptstyle}
	{\myffrac{#1}{#2} in \scriptstyle}
	{\myffrac{#1}{#2} in \scriptscriptstyle}
	{\myffrac{#1}{#2} in \scriptscriptstyle}
}

\newcommand{\wt}{\widetilde}

\newcommand{\ssup}[1] {{\scriptscriptstyle{({#1}})}}
\newcommand{\hull}{\text{\rm hull\,}}

\newcommand{\distH}{{\text{\rm dist}}_{\text{\rm H}}}

\newcommand{\frakd}{\fraktura d}

\renewcommand{\caption}[1]{%
\vglue0.5cm
\refstepcounter{figure}
\begin{center}
\begin{minipage}[c]{0.8\textwidth}\small {\sc Fig.~\thefigure\ }#1\end{minipage}
\end{center}
}

\begin{document}

\title[Shape theorem for random walks\hfill  \version \hfill]{\large Eigenvalue vs perimeter in a shape theorem \\for self-interacting random walks}

\author[\hfill  \version \hfill Biskup and Procaccia]
{Marek~Biskup$^{1,2}$ and Eviatar B.\ Procaccia$^3$}
\thanks{\hglue-4.5mm\fontsize{9.6}{9.6}\selectfont\copyright\,\textrm{2017}\ \textrm{M.~Biskup, E.B.~Procaccia.
Reproduction, by any means, of the entire
article for non-commercial purposes is permitted without charge.\vspace{2mm}}}
\maketitle

\vspace{-4mm}
\centerline{$^1$\textit{Department of Mathematics, UCLA, Los Angeles, California, USA}}
\centerline{$^2$\textit{Center for Theoretical Study, Charles University, Prague, Czech Republic}}
\centerline{$^3$\textit{Department of Mathematics, Texas A\&M University, College Station, Texas, USA}}

\begin{abstract}
We study paths of time-length~$t$ of a continuous-time random walk on~$\mathbb Z^2$ subject to self-interaction that depends on the geometry of the walk range and a collection of random, uniformly positive and finite edge weights. The interaction enters through a Gibbs weight at inverse temperature~$\beta$; the ``energy'' is the total sum of the edge weights for edges on the outer boundary of the range. For edge weights sampled from a translation-invariant, ergodic law, we prove that the range boundary condensates around an asymptotic shape in the limit~$t\to\infty$ followed by~$\beta\to\infty$. The limit shape is a minimizer (unique, modulo translates) of the sum of the principal harmonic frequency of the domain and the perimeter with respect to the first-passage percolation norm derived from (the law of) the edge weights. A dense subset of all norms in~$\mathbb R^2$, and thus a large variety of shapes, arise from the class of weight distributions to which our proofs apply. 
\end{abstract}
\maketitle

\vglue-3mm
\section{Introduction and results}
\nopagebreak\vglue-4mm
\subsection{Motivation}
\noindent
Limit theorems for random shapes have been a topic of recurring interest in both probability and statistical mechanics. One successful line of  attack came up in the 1990s under the banner of the \emph{Wulff construction}. There one was interested in the asymptotic shape of a ``droplet'' of one equilibrium phase (e.g., the minus phase of the Ising model at inverse temperature~$\beta$) immersed in another (the plus Ising phase) subject to a restriction on the overall order parameter (the magnetization, or the total number of plus spins, in the Ising case). See, e.g., Alexander, Chayes and Chayes~\cite{acc1990}, Dobrushin, Koteck\'y and Shlosman~\cite{dks1992}, Ioffe and Schonmann~\cite{is1998} for studies in the spatial dimension $d=2$ and Cerf~\cite{cerf2000}, Bodineau~\cite{bodineau1999}, Cerf and Pisztora~\cite{cp2000} and Bodineau, Ioffe and Velenik~\cite{biv2000} for work in dimensions $d\ge3$.

In all of the above mentioned studies the asymptotic shape is determined as a solution of an \emph{isoperimetric problem}; namely, 
\begin{equation}
\label{E:1.1}
\inf\bigl\{\CP(U)\colon U\subset\BR^d \text{ open},\,|U|=1\bigr\}
\end{equation}
for $|U|$ denoting the Lebesgue measure of~$U$ and $\CP(U)$ being a perimeter functional given explicitly in terms of a model-dependent surface tension $\tau\colon\BR^d\to(0,\infty)$ by
\begin{equation}
\CP(U):=\int_{\partial U}\tau\bigl(\textbf{n}(x)\bigr)H^{d-1}(\textd x),
\end{equation}
where $\textbf{n}(x)$ denotes the unit outer normal to~$\partial U$ at point~$x$ and $H^{d-1}$ is the $(d-1)$-dimensional Hausdorff measure on~$\partial U$. 
The mode of convergence of the random shape to its deterministic limit is generally stronger in $d=2$ than in~$d\ge3$. This is due to a different type of continuity of the perimeter functional. 

Recently, an interesting example which is related to but does not quite fit under the umbrella of the Wulff construction program emerged in the work of Berestycki and Yadin~\cite{by2013}. They studied random walks of time-length~$t$ that are subject to an interaction suppressing, through a Gibbs weight at inverse temperature~$\beta$, the internal vertex boundary of the walk range (i.e., the set of sites visited by the walk). The analysis in~\cite{by2013} determined the (typical) spatial size of the range in the limit $t\to\infty$ but made no definitive conclusions on the limit shape.

In the present paper we prove the existence of a limit shape for a class of self-interacting random walks closely related to that studied in~\cite{by2013}; namely, those where the ``energy'' of a random walk configuration is the length of the external edge boundary of the range measured using a sample of random, non-negative, shift-ergodic weights. A novel feature here is that this limit shape is \emph{no longer} a solution of an isoperimetric problem. Instead, after a judicious rescaling of the relevant quantities, it is the minimizer of
\begin{equation}
\label{E:1.3}
\inf\bigl\{\lambda(U)+\CP(U)\colon U\subset\BR^d\text{ open}\bigr\},
\end{equation}
where $\CP(U)$ is a quantity of the form defined above and 
\begin{equation}
\label{E:1.4}
\lambda(U):=\inf\Bigl\{\Vert\nabla g\Vert_2^2\colon g\in C^\infty(\BR^d),\,\supp(g)\subset U,
\Vert g\Vert_2=1\Bigr\}
\end{equation}
denotes the principal (i.e., smallest) eigenvalue of the negative Laplacian in~$U$ with Dirichlet boundary conditions on~$\partial U$. Some attributes of the isoperimetric problem of course remain in effect; e.g., the eigenvalue part of the functional acts as a kind of soft lower bound on the volume. However, the precise shape is a consequence of a subtle interplay between both terms and, in general, the minimizers of \eqref{E:1.1} and \eqref{E:1.3} are not homothetic. A companion paper (Biskup and Procaccia~\cite{BP-analysis}) carries out a detailed study of this variational problem.

As it turns out, our results are restricted to spatial dimension $d=2$ and they require the limit of zero temperatures, $\beta\to\infty$. The restriction on the dimension is (partly) a matter of convenience --- indeed, the control of the shape in $d\ge3$ would suffer the same technical challenges as in the above Wulff construction studies. The limit of zero temperatures is more serious and is dictated by our inability to control all aspects of the problem in full generality. (In particular, the variational problem \eqref{E:1.3} may be relevant only at~$\beta=\infty$.) Notwithstanding, since we allow the edge weights to be random, we are able to demonstrate a rich class of possible limit shapes.

\subsection{The model}
Let us now move to precise definitions. Although we will ultimately restrict attention to~$d=2$ we will keep~$d$ general in the forthcoming discussion until we get to the main result. This will allow us to track the role of~$d$ in the numerical values of some important exponents. 
 
Given a finite set~$A\subset\BZ^d$, let $\hull(A)$ denote the complement in~$\Z^d$ of the set of vertices in the unique infinite connected component of~$\BZ^d\smallsetminus A$ --- in short, $\hull(A)$ is~$A$ with all of the ``holes'' filled. Consider the continuous-time simple random walk  $\{X_t\colon t\ge0\}$ on~$\BZ^d$ with uniform jump rate~$2d$ and let $P^x$ be its law for $P^x(X_0=x)=1$. Write $\ell_t(x):=\int_0^t\1_{\{X_s=x\}}\,\textd s$
for the local time at~$x$ and define
\begin{equation}
R(t):=\hull\bigl(\{x\in\BZ^d\colon\ell_t(x)>0\}\bigr).
\end{equation}
The reason for taking the hull is that we wish the interaction energy to depend only on the geometry of the outer boundary of the visited set. 

Consider now the set $\BB(\BZ^d)$ of (unoriented) edges in~$\BZ^d$ and let~$w\colon\BB(\BZ^d)\to(0,\infty)$ be a collection of edge weights. Given a set~$A\subset\BZ^d$, let~$\partial A$ denote the set of edges in~$\BB(\BZ^d)$ that have exactly one endpoint in~$A$. The ``energy'' of a finite set~$A$ is then given by the Hamiltonian
\begin{equation}
H(A):=\sum_{e\in\partial A}w(e).
\end{equation}
For~$\beta\ge0$ we consider the Gibbs measure $Q_{\beta,t}^x$ on the path space defined by
\begin{equation}
Q_{\beta,t}^x(A):=\frac1{Z(\beta,t)}E^x\bigl(\1_A \texte^{-\beta H(R(t))}\bigr),
\end{equation}
where
\begin{equation}
\label{E:1.8}
Z(\beta,t):=E^x\bigl( \texte^{-\beta H(R(t))}\bigr)
\end{equation}
with $E^x$ denoting the expectation with respect to~$P^x$.   

A canonical (and simplest) choice of the weights is $w(e):=1$. However, this will lead to a rather uninteresting limit shape (namely, a square) and so, in order to get a larger variety of possible limit shapes, we will permit the weights to be random. Some natural restriction on the law of the weights is still in order; these are the subject of:

\begin{definition}
\label{def1}
We will call a probability law~$\BP$ on~$(0,\infty)^{\BB(\BZ^d)}$ \emph{admissible} if it is stationary and ergodic with respect to translates of~$\BZ^d$ and the marginal law of $w(e)$ is compactly supported in $(0,\infty)$ for each~$e\in\BB(\BZ^d)$. 
\end{definition}

\begin{figure}[t]
\vglue0.2cm
\centerline{\includegraphics[width=\textwidth]{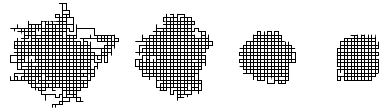}}
\begin{quote}
\small
\caption{Samples, for uniform edge weights $w(e):=1$, from $Q_{t,\beta}^0$ for the random walk of time-length $t:=5000$ at four values of the inverse temperature $\beta:=$ $0.5$, $1$, $2$ and~$5$, respectively as labeled left to right. The emergence of an asymptotic shape, which by our main result should be a perfect square in the limit $\beta\to\infty$, is quite apparent. The sizes of the external boundary in these samples are 156, 112, 66 and 54, respectively.
\label{fig1}
}
\normalsize
\end{quote}
\end{figure}

The basic problem we wish to address is the behavior of the walk sampled from $Q_{\beta,t}^0$ for~$t$ and~$\beta$ large, for a given (typical) sample of the random weights~$w$ distributed according to some admissible law~$\BP$. The main focus is the asymptotic shape of the random walk support; cf Fig.~\ref{fig1}.

A reader interested in physics might wish to have a physical system in mind that could be represented by the above class of models. One such system consists of a hydrophobic polymer chain of length~$t$ immersed in a water-based solvent. The negative affinity between the polymer and the solvent causes the polymer to fold so that the number of contacts with the solvent is minimized. The contact energy is represented by~$H(R(t))$; the inverse temperature~$\beta$ tunes the balance between the energy and the entropy of the paths. Making the weights random is quite natural, and particularly so in dimension~$d=2$, as that permits us to take into account spatial inhomogeneities of the physical substrate on which the solvent and the polymer co-exist. We refer to, e.g., the book by den Hollander~\cite{denHollander} for more information on polymers in random environment.

\subsection{Earlier and related work}
As noted above, Berestycki and Yadin~\cite{by2013} (apparently prompted by questions from I.~Benjamini) studied a related model of an interacting random walk. There are two notable differences between their and our setting: First, their interaction includes the internal components of the boundary and, second, it is given by the number of \emph{vertices} on the inner boundary. For this case they showed that the path is confined (with different type of control in $d=2$ and $d\ge3$) on the spatial scale
\begin{equation}
\label{E:1.9ue}
r(t,\beta):=\Bigl(\frac t\beta\Bigr)^{\frac1{d+1}}.
\end{equation}
The exponent is strictly less than $\ffrac12$ in all spatial dimensions $d\ge2$; the walk is thus ``squeezed'' by the interaction relative to its typical (diffusive) scaling. 

This conclusion should naturally be compared with that for the model where the ``self-inter\-action'' is proportional to the cardinality of the whole range. In this case, Bolthausen~\cite{Bolthausen} showed that the walk is confined to the spatial scale $t^{\frac1{d+2}}$. (Bolthausen's paper is actually focused on $d=2$ only but one can relate this to the problem of Brownian motion among Poissonian obstacles studied by Sznitman~\cite{S1991} in $d=2$ and Povel~\cite{P1999} in all $d\ge3$.) The minimizing shape is then determined by the variational problem
\begin{equation}
\label{E:1.10}
\inf\bigl\{\lambda(U)\colon U\subset\BR^d\text{ open},\,|U|=1\bigr\}
\end{equation}
whose unique minimizer is, thanks to the Faber-Krahn inequality, a Euclidean ball. 

The appearance of the Dirichlet eigenvalue has to do with the large-deviations cost of keeping the walk confined to a given spatial region. The associated large-deviations principle (which goes back to Donsker and Varadhan~\cite{DV1979}) underlies a large body of literature on random walks interacting via their local time and/or through an underlying random environment; e.g., the study of the parabolic Anderson model (cf K\"onig~\cite{Koenig-review} for a recent review), random walk and/or Brownian motion among random obstacles (Sznitman~\cite{Sznitman-book}), etc. Two recent papers of Asselah and Shapira~\cite{AS1,AS2} are relevant for our context as they develop a detailed large-deviation approach to the size of the boundary of the random walk range in spatial dimensions $d\ge3$. This expands on the study of moderate deviations for the Brownian ``sausage'' by van den Berg, Bolthausen and den Hollander~\cite{BBH2001}.

Our variational problem \eqref{E:1.3} is kind of a mix of \eqref{E:1.1} and \eqref{E:1.10}. This is because the mechanisms underlying the latter two problems meet, and are able to compete, in our model. Understanding this competition is key for a proper description of the walk at large time scales.

\subsection{Main results}
We are now ready to give the statements of our results. For the remainder of this paper, we will  restrict ourselves to dimension $d=2$.
The main take-away message from this note is:

\begin{theorem}[Shape theorem]
\label{thm-0}
Let~$d=2$ and suppose that the law of edge weights~$\BP$ is admissible. Then there is a non-empty bounded open convex set~$U_0\subset\BR^2$ and, for each~$\epsilon>0$, there is $\beta_0(\epsilon)<\infty$ such that for all $\beta>\beta_0(\epsilon)$ and for $r(t,\beta)$ as in \eqref{E:1.9ue},
\begin{equation}
\lim_{t\to\infty}\,
Q_{\beta,t}^0\biggl(\,\inf_{x\in\BR^d}\,\distH\Bigl(r(t,\beta)^{-1} R(t),\,x+U_0\Bigr)>\epsilon\biggr)=0
\end{equation}
for $\BP$-almost every realization of the weights.
\end{theorem}

A natural problem next is the determination of the shape~$U_0$ for specific choices of the interaction. As hinted above, this will be done by characterizing~$U_0$ as the minimizer (unique, modulo translates) of a suitable variational functional. We will take
\begin{equation}
\label{E:1.11}
\CU:=\left\{U\subset\BR^2\colon\text{Jordan domain with }\,0\in U\right\}
\end{equation}
as our class of admissible continuum domains. (This means that each $U\in\CU$ is a bounded, open, connected and simply-connected set in~$\BR^2$ whose boundary is the trace of a simple closed curve.)
Given a norm~$\rho$ on~$\BR^2$ and writing $\gamma\colon[0,1]\to\BR^2$ for the curve constituting $\partial U$ (where, necessarily, $\gamma(0)=\gamma(1)$) we define the $\rho$-perimeter of~$U$ by
\begin{equation}
\label{E:1.12}
\mathfs{P}(U):=\sup_{n\ge1}\,\,\sup_{\begin{subarray}{c}
t_0,\dots,t_n\in[0,1]\\0= t_0<t_1<\ldots<t_n=1
\end{subarray}}
\,\sum_{i=1}^n\rho\bigl(\gamma(t_i)-\gamma(t_{i-1})\bigr).
\end{equation}
For $U\in\CU$ we then set
\begin{equation}
\label{E:1.13}
\mathfs{F}(U):=\lambda(U)+\mathfs{P}(U),
\end{equation}
with $\lambda(U)$ as defined in \eqref{E:1.3}.

The object~$\CF$ plays the role of a large-deviation functional for the asymptotic shape; we will thus particularly be interested in sets in~$\CU$ that minimize~$U\mapsto\CF(U)$. This is an analysis problem that has been studied in a companion paper (Biskup and Procaccia~\cite{BP-analysis}). Let us recall the salient conclusions of this work:
\begin{enumerate}
\item[(1)] For any choice of the norm~$\rho$, the functional $U\mapsto\CF(U)$ achieves its minimum on~$\CU$. Moreover, there is $U_0\in\CU$ such that
\begin{equation}
\label{E:1.14}
\CM:=\bigl\{U\in\CU\colon \CF(U)=\min\CF\bigr\}=\bigl\{x+U_0\colon -x\in U_0\bigr\}.
\end{equation}
\item[(2)] The set $U_0$ is convex and can be taken to be symmetric under reflections ($x\mapsto -x$). It is also symmetric under all rotations that preserve the norm~$\rho$.
\item[(3)]
Letting $\distH(A,B)$ denote the Hausdorff distance between sets~$A,B\subset\BR^2$, 
for each $\epsilon>0$ there is $\delta>0$ such that, for each $U\in\CU$,
\begin{equation}
\label{E:1.15}
\distH(U,\CM)>\epsilon\quad\Rightarrow\quad \CF(U)\ge \min\CF+\delta.
\end{equation}
Moreover, the map $\rho\mapsto U_0$ is continuous as a map between the set of continuous functions on the unit circle (endowed with the supremum norm) and the set of bounded non-empty open subsets of~$\BR^2$ (endowed with the Hausdorff metric).
\end{enumerate}
In order to use these results to determine~$U_0$, we thus have to show how to extract the relevant norm~$\rho$ from the particular problem at hand. 

\smallskip
We first need some terminology and notation. Consider~$\Z^2$ regarded as a graph with its nearest-neighbor structure and let $\Z^{2\star}$ denote its graph dual. A path of length~$n$ is a sequence of vertices $x_0,\dots,x_n\in \Z^{2}$ such that $x_i$ and~$x_{i+1}$ are nearest neighbors for each~$i=0,\dots,n-1$. A similar definition applies to paths on~$\Z^{2\star}$ which we then call dual paths. A path (or a dual path) is self-avoiding if no vertex appears in the corresponding sequence of vertices more than once. 
For any $x,y\in\Z^2$, let~$\Gamma(x,y)$ denote the set of all self-avoiding dual paths $\gamma=(x_0^\star,\dots,x_n^\star)$ (for any~$n\ge0$) such that~$x_0^\star=x+(\ffrac12,\ffrac12)$ and $x_n^\star=y+(\ffrac12,\ffrac12)$.

Given a self-avoiding dual path~$\gamma=(x_0^\star,\dots,x_n^\star)$, let $(e_1,\dots,e_n)$ denote the sequence of edges in~$\Z^2$ such that~$e_i$ is the edge dual to~$(x_{i-1}^\star,x_i^\star)$. Assuming a collection of weights~$w\colon\BB(\BZ^2)\to(0,\infty)$ is given, we then set
\begin{equation}
\frakd(\gamma):=\sum_{i=1}^n w(e_i)
\end{equation}
and define
\begin{equation}
\label{E:1.18ue}
D(x,y):=\inf_{\gamma\in\Gamma(x,y)}\frakd(\gamma).
\end{equation} 
Relying on duality,~$D$ can be identified with the distance in the \emph{first-passage percolation} (FPP) for weights distributed according to~$\BP$. (In this identification, $w(e)$ is the passage time for the edge dual to~$e$.) It is known (see, e.g., Auffinger, Damron and Hanson~\cite{ADH-review} for a recent review) that the FPP distance is asymptotic to a deterministic norm. This underlies our next claim:

\begin{theorem}[Extracting the norm]
\label{thm-norm}
Suppose~$d=2$ and assume that~$\BP$ is admissible. Then there is a norm~$\rho$ on~$\BR^2$ such that for $\BP$-almost every realization of the weights,
\begin{equation}
\label{E:1.26q}
\lim_{n\to\infty}\max_{\begin{subarray}{c}
x,y\in\BZ^2\\|x|,|y|\le n
\end{subarray}}
\frac{|D(x,y) -\rho(x-y)|}n=0.
\end{equation}
\end{theorem}

The FPP norm~$\rho$ is actually constructed from --- and  is thus completely determined by --- the $\BP$-almost sure limits
\begin{equation}
\label{E:rho-def}
\lim_{n\to\infty}\frac{D(0,nx)}n = \rho(x),\qquad x\in\BZ^2,
\end{equation}
whose existence is, in our case, fairly immediate from the Subadditive Ergodic Theorem; cf Lemma~\ref{lemma-3.1ue}. 
The above statement just requires combining \eqref{E:rho-def} with some additional uniformity considerations.

As soon as the norm~$\rho$ is identified for each given law of the weights, we can complete the above shape theorem by:

\begin{theorem}[Shape characterization]
\label{thm-1}
The set~$U_0$ in Theorem~\ref{thm-0} is a minimizer (unique, up to translates) of~$U\mapsto\CF(U)$ with the perimeter functional $\CP(U)$ constructed using the asymptotic FPP norm~$\rho$ defined, e.g., by \eqref{E:rho-def}.
\end{theorem}


\subsection{Remarks and extensions}
Let us proceed with some remarks on the above setting and results. We begin with comments on the connection between the norm and the limit shape.

(1) In the canonical situation $w(e):=1$ (i.e., with constant weights),~$\rho$ is the $\ell^1$-norm and~$U_0$ is then an $\ell^\infty$-ball (a square). In the situation treated by Berestycki and Yadin~\cite{by2013} --- where the interaction is through the number of boundary vertices rather then edges --- we expect~$\rho$ to be the $\ell^\infty$-norm and~$U_0$ thus an $\ell^1$-ball (a diamond). This is in fact quite consistent with the numerical simulations shown in~\cite{by2013}.

(2) For non-degenerate laws of the weights, the FPP norm~$\rho$ is generally not explicitly computable. A question then arises what possible norms one can get in our class of models and what shapes then arise as minimizers of the associated functional $U\mapsto\CF(U)$? Here we recall: 

\begin{theorem}[H\"aggstr\"om and Meester~\cite{HM95}, Theorem 1.3 and remarks afterwards]
For each compact convex set~$B\subset\BR^d$ that has a non-empty interior and is symmetric under the reflection $x\mapsto -x$ there is a translation-invariant, ergodic law~$\BP$ on $(0,\infty)^{\BB(\BZ^d)}$ such that~$B$ is the unit ball in the norm~$\rho$ defined implicitly by \eqref{E:rho-def}. Moreover, the weights $w(e)$ are bounded under~$\BP$ and their law is strongly mixing.
\end{theorem}

The norms on~$\BR^d$ are in one-to-one correspondence with the symmetric (non-empty) compact convex sets $B\subset\BR^d$ via the equivalence
\begin{equation}
\label{E:1.21u}
B:=\{x\in\BR^d\colon\rho(x)\le1\}
\qquad\Leftrightarrow\qquad
\rho(x)=\sup\{x\cdot y\colon y\in B\},\quad\forall x\in\BR^d.
\end{equation}
In particular, all norms in~$\BR^d$ arise as FPP norms for some law of the weights as described in the theorem. However, the construction in \cite{HM95} does not make $w(e)$ bounded away from zero (the whole argument is based on existence of arbitrarily fast edges) and so we are unfortunately not guaranteed to get all norms from laws that are admissible in the sense of Definition~\ref{def1}. Notwithstanding, thanks to a continuity of the FPP norm in the underlying law of the edge weights, we get at least a dense subset of all norms.

The question on what limit shapes arise from the class of admissible laws is more complicated. Indeed, at this time we do not even know whether the map $\rho\mapsto U_0$ (of a norm to a symmetric, bounded, open and convex set) is injective and/or surjective. This is very different for the isoperimetric problem where we have the one-to-one correspondence~\eqref{E:1.21u}.

(3) A particularly interesting sub-class of admissible laws are those that make the weights i.i.d. A good amount of research has gone into the analytic properties of such norms. Durrett and Liggett~\cite{DL1981} gave examples of i.i.d.\ laws for which~$\rho$ is not strictly convex; this implies facets on the ball in the~$\rho$-metric and, by Theorem~1.7 of Biskup and Procaccia~\cite{BP-analysis}, appearance of corners on~$U_0$. However, the opposite regime --- namely, the conditions for strict convexity of the FPP norm --- is only poorly understood (for the i.i.d.\ setting); see the review by Auffinger, Damron and Hanson~\cite{ADH-review}. We advance some of the existing conjectures into:

\begin{conjecture}
\label{conj-1.5}
For i.i.d.\ laws of the edge weights that are continuously distributed and bounded uniformly away from zero and infinity, $U_0$ as well as the isoperimetric minimizer (a.k.a.\ the Wulff shape) for the associated FPP norm have no corners.  Notwithstanding, facets on both~$U_0$ and the Wulff shape are allowed and, most likely, typical.
\end{conjecture}

\noindent
The rest of our comments are concerned with possible extensions of the present work.

\smallskip
(4) We believe that a shape theorem holds for all~$\beta$ sufficiently large although we do not think that the limit shape is obtained by optimizing the functional of the form~\eqref{E:1.13}. Our need for taking $\beta\to\infty$  stems from the inability to control the surface-order terms in the ``eigenvalue part'' of the contribution to the partition function \eqref{E:1.8}; see also Remark~\ref{remark:note}. Any progress on this question is certainly of interest.

(5) Theorem~\ref{thm-0} should extend without significant changes to discrete-time random walks. Our focus on continuous-time walks is more advantageous technically as that makes functional-analytic techniques and particularly, spectral calculus, readily available. 

(6) The above results should be extendable to all $d\ge2$ (still under~$\beta\to\infty$) provided that we weaken the Hausdorff metric to
\begin{equation}
\text{dist}_1(A,B)=|A\triangle B|.
\end{equation}
The extension would carry a significant amount of technical overhead that would ultimately detract from the principal message of our paper. However, we still think it is worthy an attempt. Another obstruction is that our understanding of the variational problem for the functional \eqref{E:1.13} is presently considerably more advanced in $d=2$ than in~$d\ge3$.

(7) Theorem~\ref{thm-0} discusses the asymptotic shape of the random-walk support but one can also ask about how this shape is positioned relative to the starting point of the walk. Similarly, one might also be interested in the asymptotic law of the endpoint of the path. We expect that both of these laws can be expressed using the principal eigenfunction of the Laplacian. (Essentially, the eigenfunction should correspond to the local-time profile of the walk.)

(8) Our choice of the interaction --- the sum of non-negative weights on the outer boundary of the range --- is reasonably canonical, but other natural choices exist as well. One is the size of the inner-vertex boundary (i.e., the case related to that in Berestycki and Yadin~\cite{by2013}). However, even this example harbors some of the complications that turn out to be the more severe the more general the interaction becomes: the minimal ``energy'' of paths from $0$ to~$nx$ may no-longer be subadditive in~$n$. There are ways to overcome this in specific cases (e.g., by way of considering random-walk bridges) but a general approach seems elusive at this point.

\subsection{Main ideas and outline}
The remainder of this paper is devoted to the proofs of the above results. These come roughly in three parts each of which is the content of one section. 
Here is the basic starting idea: Consider the collection of sets
\begin{equation}
\mathfrak S:=\bigl\{S\subset\BZ^d\colon\text{ finite, connected with $S^\cc$ connected},\,0\in S\bigr\}.
\end{equation}
For the random walk started from the origin, $R(t)$ takes values in~$\mathfrak S$. Since the interaction depends only on~$R(t)$, to control the shape it suffices to study
\begin{equation}
Q_{t,\beta}^0\bigl(R(t)=S\bigr)=\frac{ \texte^{-\beta H(S)}P^0\bigl(R(t)=S\bigr)}{
\sum_{S'\in\mathfrak S} \texte^{-\beta H(S')}P^0\bigl(R(t)=S'\bigr)}.
\end{equation}
To estimate this probability, we generally need upper bounds on the numerator and lower bounds on the denominator. The essential point is that, in both of these, the computation separates into two parts: computing $P^0(R(t)=S)$ and matching the result against $ \texte^{-\beta H(S)}$. 

The proofs are then organized as follows: Section~\ref{sec3} expresses the probability $P^0(R(t)=S)$ by means of the principal eigenvalue~$\lambda^{\ssup1}_S$ of the (negative) discrete Laplacian in~$S$ (defined precisely below). This yields an expression of the form 
\begin{equation}
\label{E:4.1}
Q_{\beta,t}\bigl(R(t)=S\bigr)= \texte^{-[\,t\lambda^{\ssup1}_S+\beta H(S)+O(|\partial S|)]}.
\end{equation}
The next step, carried out in Section~\ref{sec4}, is a coarse-graining argument. This is achieved by a resummation of all~$S$ that give rise to the same continuum object, called the \emph{skeleton}~$\cmss P(S)$, which is a polygonal domain in~$\BR^2$, with interior denoted by $\Int(P)$, that closely approximates the given set~$S$ (in the Hausdorff distance). The result is now expressed in terms of the norm~$\rho$ and the perimeter functional~$\CP$ defined above in the form that looks, roughly, as
\begin{equation}
\label{E:1.26ue}
\sum_{S\colon\cmss P(S)=\cmss P}
Q_{\beta,t}\bigl(R(t)=S\bigr) =\texte^{-[\,t\lambda(\Int(\cmss P))+\beta\CP(\cmss P)+O(|\cmss P|)]}.
\end{equation}
The reason for using skeletons is that they permit us to perform local optimization of the boundary energy while, at the same time, reduce the combinatorial complexity of possible shapes. 

The expression in the exponent on the right of \eqref{E:1.26ue} has the form of the functional~$\CF$. With some extra work, this allows us to complete the proofs of main results. This occurs in Section~\ref{sec5}; the final section (Section~\ref{sec5a}) derives some useful estimates for eigenvalues and eigenfunctions of the Laplacian in finite subsets of~$\Z^2$.


\section{Extracting the principal eigenvalue}
\label{sec3}\nopagebreak\noindent
The goal of this section is to express the probability law of~$R(t)$ under the random walk measure using the principal eigenvalue of the Laplacian. This is a classical problem that lies, in one way or another, at the heart of all the aforementioned studies by Bolthausen~\cite{Bolthausen}, Sznitman~\cite{S1991}, Povel~\cite{P1999} and others. Still, since our focus is on controlling all aspects of the problem at the \emph{surface} order, some additional caution is needed.

\subsection{Key propositions}
Let $\mathfrak S$ be as above and let $S\in\mathfrak S$. Abusing the earlier notation slightly, let
\begin{equation}
\lambda^{\ssup1}_S:=\inf\,\biggl\{\frac12\sum_{\langle x,y\rangle}\bigl|g(y)-g(x)\bigr|^2\colon \supp g\subset S,\,\sum_x\bigl|g(x)\bigr|^2=1\biggr\},
\end{equation}
where the first sum is over unordered pairs of nearest-neighbor vertices in~$\BZ^d$ and the infimum is over functions $g\colon\BZ^d\to\BR$ with the stated properties. Modulo a sign change, this is the principal eigenvalue of the discrete Laplacian
\begin{equation}
\Delta f(x):=\sum_{y\colon|y-x|=1}\bigl[f(y)-f(x)\bigr]
\end{equation}
in~$S$ with Dirichlet boundary conditions on~$\BZ^2\smallsetminus S$. We will also frequently need the notation~$\lambda^{\ssup k}_S$ for the $k$-th eigenvalue of~$-\Delta$ in~$S$,  where $k=1,\dots,|S|$, and we will write $h^{\ssup k}$ to denote the associated ($k$-th) eigenfunction normalized so that
\begin{equation}
\sum_x \bigl|h^{\ssup k}(x)\bigr|^2=1.
\end{equation}
This eigenfunction is not necessarily unique, but it can always be taken real valued.

\smallskip
Let us start with bounds on~$P^0(R(t)=S)$. 
The upper bound is quite simple:

\begin{lemma}
\label{prop-3.1w}
For any $S\in\mathfrak S$,
\begin{equation}
\label{E:3.6w}
P^0\bigl(R(t)=S\bigr)\le|S|^{3/2} \texte^{-t\lambda^{\ssup1}_S}.
\end{equation}
\end{lemma}

\begin{proofsect}{Proof}
Let $\tau_S$ denote the first exit time of the random walk $X_t$ from~$S$. The function
\begin{equation}
\label{E:3.7}
f(x,t):=P^x\bigl(\tau_S>t\bigr)
\end{equation}
solves the differential equation
\begin{equation}
\frac\partial{\partial t}f(x,t)=\Delta f(x,t),\qquad x\in S,\,t>0,
\end{equation}
with initial/boundary data
\begin{equation}
f(x,0)=\ind_S(x) \quad\text{and}\quad f(\cdot,t)=0\text{ on }S^\cc.
\end{equation}
In terms of the canonical inner product $\langle\cdot,\cdot\rangle$ in~$\ell^2(S)$, we get $f(0,t)=\langle \delta_0, \texte^{t\Delta}\ind_S\rangle$ which in the language of eigenvalues/eigenfunctions of~$\Delta$ becomes
\begin{equation}
\label{E:3.9}
P^0\bigl(\tau_S>t\bigr)=\sum_{k=1}^{|S|} \texte^{-t\lambda^{\ssup k}_S} h^{\ssup k}(0)\sum_{x\in S}h^{\ssup k}(x).
\end{equation}
Invoking the Cauchy-Schwarz inequality, the bound $\lambda^{\ssup k}_S\ge\lambda^{\ssup1}_S$ along with the fact that $h^{\ssup k}$ is normalized, and in particular obeys $|h^{\ssup k}(0)|\le1$,~yield
\begin{equation}
P^0\bigl(\tau_S>t\bigr)\le |S|^{3/2} \texte^{-t\lambda^{\ssup1}_S}.
\end{equation}
But
\begin{equation}
\label{E:3.11}
P^0\bigl(R(t)=S\bigr)\le P^0\bigl(R(t)\subseteq S\bigr)= 
P^0\bigl(\tau_S>t\bigr)
\end{equation}
and so the claim follows.
\end{proofsect}

The corresponding lower bound is somewhat more involved. Fortunately, we will only need it for a reduced class of sets.  Then we have: 

\begin{proposition}
\label{prop-3.2}
Given $U\in\CU$ convex, let
\begin{equation}
\label{E:3.12}
S(U,t,\beta):=\bigl\{x\in\BZ^2\colon x/r(t,\beta)\in U\bigr\}.
\end{equation}
There is $c=c(U)\in(0,\infty)$ and, for each~$a\in(0,\infty)$, there is~$\beta_0=\beta_0(a)$ such that for all $\beta\ge\beta_0$, all $\epsilon>0$ small enough and all $S\in\mathfrak S$ with
\begin{equation}
\label{E:S-restr}
S\bigl((1-\epsilon)U,t,\beta\bigr)\subseteq S\subseteq S(U,t,\beta)\quad\text{and}\quad|\partial S|\le a\,r(t,\beta),
\end{equation}
then
\begin{equation}
\label{E:3.14w}
P^0\bigl(R(t)=S\bigr)\ge \texte^{-t\lambda^{\ssup1}_S-c|\partial S|}
\end{equation} 
holds as soon as~$t$ is sufficiently large.
\end{proposition}

Given a finite and connected~$S\subset\BZ^2$, let $\diam(S)$ denote the \emph{intrinsic diameter} of~$S$; i.e., the diameter of a graph with vertices in~$S$ and edges with both vertices in~$S$ measured using the intrinsic distance. An important input in the proof are the following estimates on the size of the spectral gap and decay of the principal eigenfunction relative to the size of the set in \eqref{E:3.12}:

\begin{lemma}
\label{lemma-2.3}
For~$U\in\CU$ convex and $\epsilon>0$ small enough, there are constants $c_1,c_2\in(0,\infty)$ such that if $S\in\mathfrak S$ obeys $S((1-\epsilon)U,t,\beta)\subseteq S\subseteq S(U,t,\beta)$, then
\begin{equation}
\label{E:2.13ue}
\lambda^{\ssup 1}_S\le \frac{c_1}{r(t,\beta)^2},
\end{equation}
\begin{equation}
\label{E:3.18a}
\lambda^{\ssup 2}_S-\lambda^{\ssup 1}_S\ge \frac{c_2}{r(t,\beta)^2}.
\end{equation}
once $r(t,\beta)$ is sufficiently large. Moreover, there are also $c_3,c_4\in(0,\infty)$ such that the principal eigenfunction in such~$S$ obeys
\begin{equation}
\label{E:3.19}
h^{\ssup1}(0)^2\ge\frac{c_3}{|S|}
\end{equation}
and
\begin{equation}
\label{E:3.17b}
\min_{z\in S} \bigl|\,h^{\ssup1}(z)\bigr|\ge \texte^{-c_4\diam(S)},
\end{equation}
as soon as $r(t,\beta)$ is sufficiently large.
\end{lemma}

The proof of this lemma requires some standard but rather technical steps that will be also used elsewhere in this paper. To avoid breaking the flow of the exposition, we therefore postpone the proof to the Appendix. The key step in the proof of Proposition~\ref{prop-3.2} is:

\begin{lemma}
\label{lemma-surface}
Let~$U\in\CU$. There is~$c=c(U)\in(0,\infty)$ and for each~$a\in(0,\infty)$ there is $\beta_0(a)\in(0,\infty)$ such that for all $\epsilon>0$ small enough, all $\beta\ge\beta_0(a)$ and any $S\in\mathfrak S$ that obeys \eqref{E:S-restr}, 
\begin{equation}
\label{E:3.18}
P^0\bigl(R(t)=S\big| R(t)\subseteq S\bigr)
\ge \texte^{-c|\partial S|}
\end{equation}
holds as soon as~$t$ is sufficiently large.
\end{lemma}

Deferring the proof of the lemma to the next subsection, let us see how it implies the above proposition:

\begin{proofsect}{Proof of Proposition~\ref{prop-3.2}}
Obviously
\begin{equation}
P^0\bigl(R(t)=S\bigr) = P^0\bigl(R(t)\subseteq S\bigr)\,P^0\bigl(R(t)=S\big| R(t)\subseteq S\bigr).
\end{equation}
In light of \eqref{E:3.9} and the equality in \eqref{E:3.11}, the first term on the right can be bounded as
\begin{equation}
\begin{aligned}
P^0\bigl(R(t)\subseteq S\bigr)
&\ge  \texte^{-t\lambda^{\ssup1}_S}h^{\ssup1}(0)^2
+\sum_{k=2}^{|S|} \texte^{-t\lambda^{\ssup k}_S} h^{\ssup k}(0)\sum_{x\in S}h^{\ssup k}(x)
\\
&\ge  \texte^{-t\lambda^{\ssup1}_S}h^{\ssup1}(0)^2-|S|^{3/2} \texte^{-t\lambda^{\ssup 2}_S},
\end{aligned}
\end{equation}
where we applied that $h^{\ssup1}$ is of one sign and then bounded the contribution of the second and higher eigenvalues as in the proof of Lemma~\ref{prop-3.1w}.
Invoking \twoeqref{E:3.18a}{E:3.19} along with the definition of $r(t,\beta)$, we  then see that, for some $c>0$,
\begin{equation}
P^0\bigl(R(t)\subseteq S\bigr)\ge\frac{c}{|S|} \texte^{-t\lambda^{\ssup1}_S}
\end{equation}
once~$t$ is sufficiently large. In combination with Lemma~\ref{lemma-surface}, this proves the claim.
\end{proofsect}

\subsection{Surface order in confinement probability}
Let us now move to the proof of Lemma~\ref{lemma-surface}. As it turns out, the main challenge is to deal with the conditioning on~$R(t)\subseteq S$, which we can write as $\tau_S>t$ where (we recall) $\tau_S:=\inf\{t\ge0\colon X_t\not\in S\}$. This would perhaps appear easier if the condition $\tau_S>t$ were replaced by $\tau_S=\infty$. As it turns out, this is not hard to arrange:

\begin{lemma}
\label{lemma-2.5a}
There are constants $c,\tilde c\in(0,\infty)$ depending only on~$U$ such that for all $\epsilon>0$ small, all $\beta\ge1$ and all~$a\in(0,\infty)$, if $S\in\mathfrak S$ obeys \eqref{E:S-restr}, then 
\begin{equation}
P^0(A|\tau_S>t)\ge\frac12 P^0(A|\tau_S=\infty)- \texte^{-ct}
\end{equation}
holds for all events $A\in\sigma(X_u\colon 0\le u\le t-u(t))$ where
\begin{equation}
\label{E:2.23}
u(t):=\tilde c\, a\, r(t,\beta)^3.
\end{equation}
\end{lemma}

\begin{proofsect}{Proof}
Let~$u(t)$ be as above with a constant~$\tilde c$ to be determined.
Abbreviate $t':=t-u(t)$ and let $A\in\sigma(X_u\colon0\le u\le t')$. Define the collection of stopping times $\{T_i\colon i\ge0\}$ inductively by
\begin{equation}
T_{i+1}:=t'\wedge\inf\{u>T_i\colon X_u\ne X_{T_i}\}, \quad\text{where}\quad T_0:=0,
\end{equation}
and denote 
\begin{equation}
N:=\inf\{i\in\BN\colon T_i=t'\}.
\end{equation}
Using the function $f(x,t)$ from \eqref{E:3.7}, the Markov property gives us
\begin{equation}
P^0(A|\tau_S>t) = 
E^0\biggl(\1_{A\cap\{\tau_S>t'\}}\prod_{i=1}^N\frac{f(X_{T_i},t-T_i)}{f(X_{T_{i-1}},t-T_{i-1})}\biggr).
\end{equation}
The function $f(x,t)$ admits the representation \eqref{E:3.9}. The estimates we derived for the $k\ge2$ part of the sum above show
\begin{equation}
\biggl|\,f(x,t)-\texte^{-t\lambda^{\ssup1}_S}h^{\ssup1}(x)\sum_{z\in S}h^{\ssup1}(z)\biggr|
\le |S|^{3/2}\texte^{-t\lambda^{\ssup2}_S}.
\end{equation}
For $0\le u\le t$ we then get
\begin{equation}
\frac{f(x,t-u)}{f(y,t)}\ge \texte^{u\lambda^{\ssup1}_S}\frac{h^{\ssup1}(x)}{h^{\ssup1}(y)}\biggl(1-2\frac{|S|^{3/2}}{\min_{z\in S} h^{\ssup1}(z)^2} \texte^{-(t-u)(\lambda^{\ssup2}_S-\lambda^{\ssup1}_S)}\biggr).
\end{equation}
Thanks to \eqref{E:3.18a} and \eqref{E:3.17b} we can bound
\begin{equation}
\label{E:2.29ua}
\frac{|S|^{3/2}}{\min_{z\in S} h^{\ssup1}(z)^2} \texte^{-(t-u)(\lambda^{\ssup2}_S-\lambda^{\ssup1}_S)}
\le |S|^{3/2}\texte^{2c_4 \diam(S)-u(t)r(t,\beta)^{-2}},\qquad 0\le u\le t'.
\end{equation}
Using the bound on~$|\partial S|$ along with the fact that
\begin{equation}
|\partial S|\ge\diam(S)
\end{equation}
and also the fact that $|S|=O(r(t,\beta)^2)$, the right-hand side of \eqref{E:2.29ua} decays to zero with~$t\to\infty$ as soon as $\tilde c$ is take sufficiently large. We will henceforth assume~$t$ is large enough so that the expression on the right of \eqref{E:2.29ua} is less than~$\ffrac14$. 

For all paths of the walk in the event $\{N\le2t\}$, we thus get
\begin{equation}
\prod_{i=1}^N\frac{f(X_{T_i},t-T_i)}{f(X_{T_{i-1}},t-T_{i-1})}
\ge\frac12 \texte^{t'\lambda^{\ssup1}_S}\frac{h^{\ssup1}(X_{t'})}{h^{\ssup1}(0)},
\end{equation}
and thus
\begin{equation}
P^0(A|\tau_S>t) \ge\frac12  \texte^{t'\lambda^{\ssup1}_S}\frac1{h^{\ssup1}(0)}E^0\bigl(\1_{A\cap\{\tau_S>t'\}\cap\{N\le 2t\}} h^{\ssup1}(X_{t'})\bigr).
\end{equation}
But the Doob $h$-transform (or a limit argument based on the spectral decomposition) gives 
\begin{equation}
\label{E:3.28}
P^0(A|\tau_S=\infty) =  \texte^{t'\lambda^{\ssup1}_S}\frac1{h^{\ssup1}(0)}E^0\bigl(\1_{A\cap\{\tau_S>t'\}} h^{\ssup1}(X_{t'})\bigr)
\end{equation}
so we just need to prove
\begin{equation}
P^0\bigl(N>2t\big|\tau_S=\infty\bigr)\le \texte^{-ct}.
\end{equation}
This is derived directly from \eqref{E:3.28}, \eqref{E:3.19} and the fact that $t\lambda^{\ssup1}_S=o(t)$; cf~\eqref{E:2.13ue}.
\end{proofsect}

In order to proceed further, it will be easier to convert to a discrete-time version of the conditional chain. Writing temporarily~$d$ for the spatial dimension, here we note:

\begin{lemma}
\label{lemma-2.6a}
Under $P^0(-|\tau_S=\infty)$, the law of $(X_t\colon t\ge0)$ satisfes
\begin{equation}
(X_t\colon t\ge0)\,\,\overset{\text{\rm law}}=\,\, (Z_{N_t}\colon t\ge0),
\end{equation}
where $Z=(Z_n\colon n\ge0)$ is the discrete time Markov chain on~$S$ with $Z_0=0$ a.s.\ and the transition probabilities given by
\begin{equation}
\label{E:2.29ue}
\cmss P(x,y)=\frac1{2d-\lambda^{\ssup1}_S}\,\,\frac{h^{\ssup1}(y)}{h^{\ssup1}(x)}
\end{equation}
whenever $x$ and~$y$ are nearest neighbors of~$\BZ^d$, and $(N_t\colon t\ge0)$ is the Poisson process with rate $2d-\lambda^{\ssup1}_S$, independent of~$Z$.
\end{lemma}

\begin{proofsect}{Proof}
The general theory ensures that the random walk conditioned on $\tau_S=\infty$ is still a Markov process. Since the state space is finite, this Markov process is automatically realized as a time change of a discrete-time Markov chain. An inspection of \eqref{E:3.28} shows that the waiting time at site~$x$ is exponential with parameter $2d-\lambda^{\ssup1}_S$ which, in particular, is independent of~$x$. The same formula shows that the corresponding discrete-time Markov chain has transition probabilities as stated in~\eqref{E:2.29ue}.
\end{proofsect}

Henceforth, we will think of~$X$ as realized by~$Z$ and~$N$. In this representation, there is no need to impose the conditioning on $\tau_S=\infty$ as the Markov chain~$Z$ never leaves~$S$ with probability one. Before we get to the proof of Lemma~\ref{lemma-surface}, we need one more observation:

\begin{lemma}
\label{lemma-2.7}
For any~$S\in\mathfrak S$ and any pair of nearest neighbors~$x$ and~$y$ in~$S$,
\begin{equation}
\frac{h^{\ssup1}(y)}{h^{\ssup1}(x)}\ge\frac1{2d}.
\end{equation}
\end{lemma}

\begin{proofsect}{Proof}
Since~$S$ is connected, we may assume that $h^{\ssup1}$ is of one sign, say, $h^{\ssup1}\ge0$. Whenever $x$ and~$y$ are neighbors, we thus get
\begin{equation}
h^{\ssup1}(x)\le\sum_{z\colon|z-y|=1}h^{\ssup1}(z)=(2d-\lambda^{\ssup1}_S)h^{\ssup1}(y)\le 2d\,h^{\ssup1}(y),
\end{equation}
using the eigenvalue equation and the fact that $\lambda^{\ssup1}_S\ge0$.
\end{proofsect}

\begin{proofsect}{Proof of Lemma~\ref{lemma-surface}}
Set $\beta_0(a):=\max\{1,2\tilde c a\}$, where $\tilde c$ is as in Lemma~\ref{lemma-2.5a}. In light of \eqref{E:1.9ue}, the quantity in \eqref{E:2.23} obeys  $u(t)\le t/2$ as soon as $\beta\ge\beta_0(a)$. In particular, as soon as $S$ obeys \eqref{E:S-restr} we have
\begin{equation}
P^0\bigl(R(t)=S\big| R(t)\subseteq S\bigr)
\ge \frac12 P^0\bigl(R(t/2)=S\big|\tau_S=\infty\bigr)-\texte^{-ct}.
\end{equation}
The conditional~$X$ process is now representable using~$Z$ and~$N$ above. Since the latter processes are independent, dropping probability of order $ \texte^{-ct}$ for some $c>0$ small, we can assume that the Markov chain~$X$ makes at least~$n:=ct$ steps in time~$t/2$. (This uses that the jump rate $2d-\lambda^{\ssup1}_S$ is uniformly positive thanks to \eqref{E:2.13ue}.) Hence,
\begin{equation}
P^0\bigl(R(t/2)=S\big|\tau_S=\infty\bigr)\ge P^0\bigl(\partial^{\star} S\subseteq\{Z_1,\dots,Z_n\}\bigr)-\texte^{-ct},
\end{equation}
where $\partial^{\star} S$ denotes the inner vertex boundary of~$S$.
Noting that, under our restrictions in \eqref{E:S-restr}, both $\diam(S)$ and~$|\partial^{\star} S|$ are at most order~$t^{1/3}$ which is much less than~$n$, we now realize the lower bound by forcing~$Z$ to take a path going from~$0$ directly to $\partial^{\star} S$ and then move on or around~$\partial^{\star} S$ until all vertices of~$\partial^{\star} S$ are visited. The length of this path is at most a constant times~$|\partial S|$. By Lemma~\ref{lemma-2.7}, $\cmss P(x,y)\ge c>0$ uniformly for all neighbors~$x,y\in S$, and so the probability of the given strategy is at least~$ \texte^{-c|\partial S|}$. By \eqref{E:S-restr} again, this is sufficiently large to absorb all prefactors as well as $\texte^{-ct}$ corrections.
\end{proofsect}


\begin{remark}
\label{remark:note}
The boundary term $c|\partial S|$ in the exponent on the right-hand side of \eqref{E:3.14w} arises directly from \eqref{E:3.18} (other errors are subexponential in the boundary); no such term appears in the upper bound in \eqref{E:3.6w}. Our inability to control this term more explicitly than by way of an estimate is the \emph{sole reason} why our proof of the shape theorem is restricted to the limit~$\beta\to\infty$. 

An interested reader might wonder (as we have) whether the boundary-order error term is not just an artifact of our way of proving the lower bound. The answer is that, most likely, it is not: Indeed, denoting, for $\epsilon>0$,
\begin{equation}
\mathfrak S_{t,\epsilon}':=\Bigl\{S\in\mathfrak S\colon P^0\bigl(R(t)=S\big| R(t)\subseteq S\bigr)
\ge \texte^{-\epsilon|\partial S|}\Bigr\},
\end{equation}
we have
\begin{equation}
Z(\beta,t)\ge\sum_{S\in\mathfrak S_{t,\epsilon}'} \texte^{-\beta H(S)-(\epsilon+o(1))|\partial S|-t\lambda_S^{\ssup 1}},
\end{equation}
where $o(1)$ tends to zero with diameter of~$S$. However, $H(S)=O(|\partial S|)$ and so, for~$\beta\ll\epsilon\ll1$, as soon as $\mathfrak S_{t,\epsilon}'$ contains, say, all sets in~$\mathfrak S$ within $\log r(t,\beta)$ Hausdorff distance of the Euclidean ball of radius $Kr(t,\beta)$ for some~$K>0$ large enough, the sum would diverge in the limit as~$t\to\infty$. This would be in contradiction with $Z(\beta,t)\le1$, implied by~$H(S)\ge0$.
\end{remark}

\section{Boundary norm and skeletons}
\label{sec4}\nopagebreak\noindent
Our next task is to control the energy of each configuration using a suitable norm. In addition, we will group  the  discrete ``shapes'' into (generally non-disjoint) families each of which is represented by a single continuum object, called a skeleton. A similar coarse-graining argument has been invoked numerous times in various studies of the two-dimensional Wulff construction. We will draw on ideas from the early work in the context of the two-dimensional Ising model at low temperatures by Dobrushin, Koteck\'y and Shlosman~\cite{dks1992}. 

\subsection{Boundary norm}
\noindent
The starting point is the construction of the norm~$\rho$ for each admissible law~$\BP$. Here and henceforth, we will write $\BE$ to denote the expectation with respect to~$\BP$. Recall the definition of $D(x,y)$ from \eqref{E:1.18ue}. The following lemma is standard (see, e.g., Kesten~\cite{kesten1986aspects}); we include it for completeness of exposition.

\begin{lemma}
\label{lemma-3.1ue}
Suppose~$\BP$ is an admissible law and let $0<a<b<\infty$ be such that the edge weights are supported in~$[a,b]$ under~$\BP$. Then there exists a norm~$\rho$ on~$\BR^2$ such that
\begin{equation}
\label{E:3.1ue}
\lim_{n\to\infty}\frac{\BE \,D(0,nx)}n=\rho(x),\qquad x\in\BZ^2.
\end{equation}
In particular, the limit exists and the norm~$\rho$ obeys $a|x|_1\le\rho(x)\le b|x|_1$ for all~$x$, where~$|\cdot|_1$ denotes the $\ell^1$-norm on~$\BZ^2$. 
\end{lemma}

\begin{proofsect}{Proof}
By definition,
\begin{equation}
\label{E:3.1u}
D(x,y)\le D(x,z)+D(z,y),\qquad x,y,z\in\BZ^2.
\end{equation}
The stationarity of~$\BP$ with respect to translates of~$\BZ^2$ ensures that $D(x,x+y)$ is equidistributed to~$D(0,y)$. Moreover, we clearly have
\begin{equation}
\label{E:3.2u}
a|x-y|_1\le D(x,y)\le b|x-y|_1,
\end{equation}
It follows that, for each~$x\in\BZ^2$, the sequence $\{\BE D(0,nx)\colon n\ge1\}$ is subadditive. By Fekete's Lemma the limit in \eqref{E:3.1ue} exists and defines a map $\rho\colon\Z^2\to[0,\infty)$. By \eqref{E:3.1u},~$\rho$ obeys the triangle inequality on~$\BZ^2$ and, since it is also homogeneous over integers, it is a restriction of a norm on~$\BR^2$. The bounds between~$\rho$ and the $\ell^1$-norm then follow from \eqref{E:3.2u}.
\end{proofsect}

We are now ready to address the main part of Theorem~\ref{thm-norm} which states, roughly, that the norm~$\rho$ captures the leading-order behavior of the FPP-distances $D(x,y)$ at linear scales in sufficiently large boxes centered at the origin.

\begin{proofsect}{Proof of Theorem~\ref{thm-norm}}
The inequality \eqref{E:3.1u} and the fact that $\{D(z+x,z+y)\colon z\in\BZ^2\}$ is stationary (and bounded)  implies via the Subadditive Ergodic Theorem that the limit $\hat\rho(x):=\lim_{n\to\infty}\frac{D(0,nx)}n$ exists $\BP$-almost surely. Thanks to the triangle inequality \eqref{E:3.1u}, one checks that $\hat\rho$ is invariant (a.s.) under translations of~$\BZ^2$. By ergodicity of~$\BP$, it is thus almost surely constant. Taking expectations, the Bounded Convergence Theorem and Lemma~\ref{lemma-3.1ue} above show that
\begin{equation}
\label{E:3.4u}
\lim_{n\to\infty}\frac{D(0,nx)}n = \lim_{n\to\infty}\frac{\BE \,D(0,nx)}n =\rho(x)
\end{equation}
holds $\BP$-almost surely for each~$x\in\BZ^2$. What remains to be done is to bootstrap this ``directional'' convergence into the uniform bound \eqref{E:1.26q}.

First note that, by choosing a finite set of rational directions on the unit circle and using \eqref{E:3.1u} to interpolate the directions in-between, we can augment \eqref{E:3.4u} to the following claim: For each $x\in\BZ^2$ and each~$\epsilon>0$, there is $n_0(x,\epsilon)$ with $\BP(n_0(x,\epsilon)<\infty)=1$ such that for each~$y\in\BZ^2$,
\begin{equation}
\label{E:3.5}
|y-x|_1>n_0(x,\epsilon)\quad\Rightarrow\quad (1-\epsilon)\rho(x-y)\le D(x,y)\le (1+\epsilon)\rho(x-y).
\end{equation}
(This is standard as it is exactly what one needs to conclude a uniform shape theorem for large balls in first passage percolation.) The sequence $\{n_0(x,\epsilon)\colon x\in\BZ^2\}$ is stationary and, since it is also a.s.\ finite, there is an $M>0$ such that $\BP(n_0(x,\epsilon)\le M)>0$ for all~$x\in\BZ^2$. The ergodicity of~$\BP$ then ensures that
\begin{equation}
A_M:=\{x\in\Z^2\colon\,n_0(x,\epsilon)\le M\}
\end{equation}
has a positive density in~$\BZ^2$. This implies 
\begin{equation}
\lim_{n\to\infty}
\,\frac1n\,\max_{z\colon|z|_1\le n}\dist_1(\,z,A_M)=0
\end{equation}
$\BP$-almost surely, where $\dist_1$ is the distance in the $\ell^1$-norm. 
We will now fix $\delta>0$, assume that~$n$ is so large that
\begin{equation}
\label{E:3.8u}
\delta n\ge M\quad\text{and}\quad
\max_{z\colon|z|_1\le n}\dist_1(\,z,A_M)\le\delta n
\end{equation}
holds and derive a bound on the maximum in the statement of the theorem.

Suppose that~\eqref{E:3.8u} holds and pick $x,y\in\BZ^2$ with $|x|_1,|y|_1\le n$. Let $x_0$, resp., $y_0$ denote the closest point in $\ell^1$-distance in~$A_M$ to~$x$, resp.,~$y$. The triangle inequality \eqref{E:3.1u} shows
\begin{equation}
\bigl|\,D(x,y)-D(x_0,y_0)\bigr|\le D(x,x_0)+D(y,y_0).
\end{equation}
Letting $b$ denote the a.s.\ upper bound on the weights under~$\BP$, we have
\begin{equation}
|x-x_0|_1\le M\quad\Rightarrow\quad D(x,x_0)\le b|x-x_0|_1\le bM
\end{equation}
while, assuming $\epsilon<1$, \eqref{E:3.5} and \eqref{E:3.8u} yield
\begin{equation}
|x-x_0|_1>M\quad\Rightarrow\quad
D(x,x_0)\le (1+\epsilon)\rho(x-x_0)\le 2b|x-x_0|_1\le 2b\delta n.
\end{equation}
Completely analogous statements hold for $D(y,y_0)$ as well. Using $\delta n\ge M$ we thus get
\begin{equation}
\label{E:3.12u}
\text{\eqref{E:3.8u} holds}\quad\Rightarrow\quad \bigl|\,D(x,y)-D(x_0,y_0)\bigr|\le 4b\delta n.
\end{equation}
It remains to relate $D(x_0,y_0)$ to~$\rho(x-y)$. Invoking again \eqref{E:3.5}, a very similar reasoning to one above yields
\begin{equation}
(1-\epsilon)\rho(x_0,y_0)-b\delta n\le D(x_0,y_0)\le (1+\epsilon)\rho(x_0,y_0)+b\delta n.
\end{equation}
Since \eqref{E:3.2u} gives
\begin{equation}
\label{E:3.14u}
\bigl|\rho(x-y)-\rho(x_0-y_0)\bigr|\le 2b\delta n
\end{equation}
from \twoeqref{E:3.12u}{E:3.14u} we thus conclude
 \begin{equation}
\text{\eqref{E:3.8u} holds}\quad\Rightarrow\quad \bigl|D(x,y)-\rho(x-y)\bigr|\le b(\epsilon+7\delta)n,
\quad |x|_1,|y|_1\le n.
\end{equation}
As~$\epsilon$ and~$\delta$ are arbitrary, this proves the claim.
\end{proofsect}

A slight complication with using the quantity~$D(x,y)$ is that, in our arguments, the limit~$\beta\to\infty$ has to be taken only after~$t\to\infty$. This requires us to work with the model at a fixed finite (albeit arbitrarily large)~$\beta$. The following lemma will be quite useful:

\begin{lemma}
\label{lemma-4.3}
Suppose~$\BP$ is admissible and let $\kappa>0$ be an almost-sure lower bound on the weights. Pick $\beta_0>\frac1{\kappa}\log4$. Then for all $\beta>\beta_0$ and all~$x,y\in\Z^2$,
\begin{equation}
 \texte^{-\beta D(x,y)}\le
\sum_{\gamma\in\Gamma(x,y)} \texte^{-\beta\frakd(\gamma)}\le \bigl(1-4\texte^{-\beta_0 \kappa}\bigr)^{-1}\,\texte^{-(\beta-\beta_0) D(x,y)}.
\end{equation}
\end{lemma}

\begin{proofsect}{Proof}
The lower bound is immediate (and we state it mostly for \ae sthetic reasons). For the upper bound, we first dominate
\begin{equation}
\texte^{-\beta\frakd(\gamma)}\le\texte^{-(\beta-\beta_0)D(x,y)}\texte^{-\beta_0\frakd(\gamma)},\qquad\gamma\in\Gamma(x,y),
\end{equation}
 and then use $\frakd(\gamma)\ge \kappa |\gamma|_1$, where $|\gamma|_1$ is the length of (i.e., the number of edges crossed by)~$\gamma$, and a standard path-counting argument to estimate the sum of~$\texte^{-\beta_0\frakd(\gamma)}$ over~$\gamma\in\Gamma(x,y)$ by a geometric series with quotient $4\texte^{-\beta_0 \kappa}$.
\end{proofsect}

\subsection{Polygonal approximation}
Consider now a set $S\in\mathfrak S$ and recall that~$\partial S$ is the set of outer-boundary edges of~$S$. The collection of edges dual to those in~$\partial S$ can then be ordered in such a way that they form a simple closed path~$\gamma_S$ on~$\BZ^{2\star}$. Our next task is to partition this path into pieces whose contribution can be represented using the norm constructed earlier.

Fix an integer~$L\ge1$. Given a set~$S\in\mathfrak S$ and the boundary curve~$\gamma_S$, we will define a collection of points $\{x_0,x_1,\dots,x_{n(S)}=x_0\}$ on~$\gamma_S$ inductively as follows: Let $x_0$ be the vertex of the $\Z^{2\star}$ that is the smallest of all vertices on~$\gamma_S$ in the standard lexicographic order on~$\Z^2$. Following~$\gamma_S$ in the counterclockwise orientation, let $x_2$ be the first vertex on~$\gamma_S$ that is at least $\ell^1$-distance~$L$ from~$x_1$; if no such vertex exists we set $n(S):=1$ and let $x_{n(S)}:=x_0$. Similarly, if $\{x_0,\dots,x_k\}$ have already been defined, we let $x_{k+1}$ be the first vertex on the portion of~$\gamma_S$ starting from~$x_k$ that is $\ell^1$-distance~$L$ from~$x_k$; if no such vertex exists we set~$n(S):=k+1$ and~$x_{n(S)}:=x_0$. 

Now assume that $S\in\mathfrak S$ is such that~$n(S)>1$. Connecting the pairs $(x_i,x_{i+1})$ by linear segments gives rise to a polygonal closed curve $\cmss P(S)$. In analogy with the objects defined in Dobrushin, Koteck\'y and Shlosman~\cite{dks1992}, we will refer to $\cmss P(S)$ as the $L$-\emph{skeleton} of~$S$ (the reference to~$L$ marks the obvious dependence of~$\cmss P(S)$ on~$L$). Abusing notation slightly, we write $n(\cmss P)$ for the number of vertices in~$\cmss P$. Let
\begin{equation}
\mathfrak P_L:=\bigl\{\cmss P(S)\colon S\in\mathfrak S,\,n(S)>1\bigr\}
\end{equation}
denote the set of all (non-trivial) $L$-skeletons that can possibly arise. The skeletons will enter the computation of the probabilities of various shapes via the following estimate: 

\begin{lemma}
\label{lemma-4.4}
For each admissible law~$\BP$ and each~$\epsilon\in(0,1)$ there is $\beta_1\in(0,\infty)$ and a (random) $L_0\in[1,\infty)$ such that for all $\beta>\beta_1$, all $L\ge L_0$ and all $\cmss P\in\mathfrak P_L$ with $n:=n(\cmss P)\le 1/\epsilon$,
\begin{equation}
\label{E:4.13}
\sum_{S\colon\cmss P(S)=\cmss P}
 \texte^{-\beta H(S)}\le  \texte^{2b\beta L}\exp\Bigl\{-(1-\epsilon)\sum_{i=1}^{n}\rho(x_i-x_{i-1})\Bigr\}
\end{equation}
where $x_1,\dots,x_{n}$ are the vertices of the polygonal curve~$\cmss P$ and where $b$ is an a.s.\ upper bound on the edge weights under~$\BP$.
\end{lemma}

\begin{proofsect}{Proof}
Fix~$\epsilon\in(0,1)$ and let $L_0$ be the smallest non-negative integer such that for all $L\ge L_0$ and all~$x,y\in\BZ^2$, we have
\begin{equation}
|x|_1,|y|_1\le L/\epsilon\quad\&\quad |x-y|\ge L\quad\Rightarrow\quad D(x,y)\ge(1-\epsilon)^{1/2}\rho(x-y).
\end{equation}
Theorem~\ref{thm-norm} ensures that $\BP(L_0<\infty)=1$. Let now~$L\ge L_0$ and let $\cmss P\in\mathfrak P_L$ obey $n:=n(\cmss P)\le 1/\epsilon$. Note that the vertices $x_0,\dots,x_n$ of~$\cmss P$ then satisfy $|x_i|_1\le L/\epsilon$ and, in particular, we have
\begin{equation}
\label{E:3.21ue}
D(x_{i-1},x_i)\ge (1-\epsilon)^{1/2}\,\rho(x_i-x_{i-1}),\qquad i=1,\dots,n-1.
\end{equation}
We now fix the skeleton~$\cmss P$ and proceed to prove \eqref{E:4.13} for~$\beta_1$ defined (later) only in terms of $\epsilon$ and the a.s.\ upper and lower bounds on the edge weights.

Given $S\in\mathfrak S$ such that $\cmss P(S)=\cmss P$, let $\gamma_S$ be the simple path on dual-$\Z^2$ corresponding to~$S$ as discussed above. Then $\gamma_S$ passes through the vertices $x_1,\dots,x_n$ of~$\cmss P$ in the given order and so we may write~$\gamma_i$ to denote the portion of~$\gamma_S$ between $x_{i-1}$ and~$x_i$. Clearly, $\gamma_i\in\Gamma(x_{i-1},x_i)$ and so we have an injection
\begin{equation}
\bigl\{S\in\mathfrak S\colon \cmss P(S)=\cmss P\bigr\}\,\,\to\,\,\Gamma(x_0,x_1)\times\dots\times\Gamma(x_{n-1},x_n).
\end{equation}
Moreover,
\begin{equation}
H(S) \ge \frakd(\gamma_1)+\dots+\frakd(\gamma_{n}).
\end{equation}
It thus follows that
\begin{equation}
\sum_{S\colon\cmss P(S)=\cmss P}
 \texte^{-\beta H(S)}
\le\prod_{i=1}^{n}\,\sum_{\gamma\in\Gamma(x_{i-1},x_i)} \texte^{-\beta\frakd(\gamma)}.
\end{equation}
Invoking the upper bound in Lemma~\ref{lemma-4.3}, for $\beta>\beta_0$ we then get
\begin{equation}
\sum_{S\colon\cmss P(S)=\cmss P}\texte^{-\beta H(S)}\le (1-4\texte^{-\beta_0a})^{-n}\,\prod_{i=1}^{n}\texte^{-(\beta-\beta_0)D(x_{i-1},x_i)}.
\end{equation}
Assuming that $\beta_1$ is so large that $1-\beta_0/\beta_1\ge(1-\epsilon)^{1/2}$, as soon as $\beta>\beta_1$ the bound \eqref{E:3.21ue} yields \eqref{E:4.13} with the right-hand side multiplied by
\begin{equation}
C_L :=  \texte^{\beta\rho(x_n-x_{n-1})-2\beta bL}(1-4\texte^{-\beta_0a})^{-n}.
\end{equation} 
In light of $n\le1/\epsilon$ and $\rho(x_n-x_{n-1})\le b|x_n-x_{n-1}|_1\le b L$, we have $C_L\le1$ as soon as~$\beta_1$ is so large that $\texte^{\beta_1 b}\ge (1-4\texte^{-\beta_0a})^{-1/\epsilon}$.
\end{proofsect}

A minor problem with the above construction of~$L$-skeletons is that~$\cmss P(S)$ need not be a simple curve. However, if we define $w_{\cmss P}(x)$ to be the winding number of~$\cmss P$ (oriented counterclockwise) around~$x\in\BR^2\smallsetminus\cmss P$ and set
\begin{equation}
\Int(\cmss P):=\bigl\{x\in\BR^2\smallsetminus\cmss P\colon w_{\cmss P}(x)\text{ odd}\bigr\},
\end{equation}
then $\cmss P$ can be reparametrized as a closed curve that ``keeps~$\Int(\cmss P)$ on its left'' --- i.e., after reparametrization, all points in $\Int(\cmss P)$ have a positive winding number. (An explicit construction requires a limit argument performed, e.g., in~\cite[Lemma~4.4]{blpr2014}.) In this case, it is quite consistent to abbreviate
\begin{equation}
\CP\bigl(\Int(\cmss P)\bigr):=\sum_{i=1}^{n}\rho(x_i-x_{i-1}),
\end{equation}
where $x_0,x_1,\dots,x_n=x_0$ are the vertices of~$\cmss P$.

Of course, in light of \eqref{E:4.1}, the bound in \eqref{E:4.13} will only be useful if we can first lower-bound~$\lambda_S^{\ssup 1}$ by a quantity that depends only on $\cmss P:=\cmss P(S)$. Naturally, we would like to work with~$\Int(\cmss P)$, but the fact that~$\partial S$ may reach as far as~$L$ lattice spacings outside of~$\Int(\cmss P)$ forces us to consider instead an $L$-neighborhood of~$\Int(\cmss P)$. Relating the perimeter of this~$L$-neighborhood to $\CP(\Int(\cmss P))$ then seems rather hard, due to possible ``holes'' and and other effects caused by ``inward'' spikes of~$\cmss P$. These issues disappear when we consider the convex hull of~$\cmss P$,
\begin{equation}
\hull(\cmss P):=\bigl\{\alpha x+(1-\alpha) y\colon x,y\in\cmss P,\,\alpha\in[0,1]\bigr\},
\end{equation}
and, given~$L\ge1$, let
\begin{equation}
\label{E:UL}
V_L=V_L(\cmss P):=\bigl\{x\in\BR^2\colon\dist_1(x,\hull(\cmss P))<L+3\bigr\},
\end{equation}
where (we recall) $\dist_1$ denotes the $\ell^1$-distance on~$\BR^2$. Obviously,
\begin{equation}
S\subset V_L\bigl(\cmss P(S)\bigr),\qquad S\in\mathfrak S.
\end{equation}
Concerning the relation between $\lambda_S$ and $\lambda(V_L)$, we get

\begin{proposition}
\label{prop-4.5}
There is a constant~$C\in(0,\infty)$ such that the following holds: For any $L\ge1$ and any $S\in\mathfrak S$ with $\cmss P(S)=\cmss P$ and $V_L=V_L(\cmss P)$ as above,
\begin{equation}
\label{E:4.19}
C\sqrt{\lambda_S^{\ssup1}}<1\quad\Rightarrow\quad
\frac{\lambda_S^{\ssup1}}{\bigl(1-C\sqrt{\lambda_S^{\ssup1}}\bigr)^2}\ge\lambda(V_L).
\end{equation}
\end{proposition}

We note that, since $\lambda_S^{\ssup1}\le c|S|^{-1/2}$ for some constant~$c\in(0,\infty)$, the conclusion of \eqref{E:4.19} applies as soon as~$|S|$ is sufficiently large. For the proof we will need:

\begin{lemma}
\label{fe}
There is a constant~$C\in(0,\infty)$ for which the following holds: For
any $f\colon\BZ^2\to\BR$ there is $\wt f\colon \BR^2\to\BR$ such that
\settowidth{\leftmargini}{(11)}
\begin{enumerate}
\item the map $f\mapsto\wt f$ is linear,
\item $\wt f$ is continuous on~$\BR^2$ and $\wt f(x)=f(x)$ for all $x\in \Z^2$,
\item
for any $x\in\Z^d$ and any $y\in x+[0,1)^2$ we have
\begin{equation}
\label{E:3.13aa}
\bigl|\wt f(y)\bigr|\le\max_{z\in x+\{0,1\}^2}\,\bigl|f(z)\bigr|,
\end{equation}
\item
the $L^2$-norms of the functions are related by
\begin{equation}
\label{E:3.16q}
\Bigl|\,\Vert \wt f\Vert_{L^2(\BR^2)}-\Vert f\Vert_{\ell^2(\Z^2)}\Bigr|\le C\|\nabla 
f\|_{\ell^2(\Z^2)},
\end{equation}
\item
$\wt f$ is piece-wise linear and thus a.e.\ differentiable with
\begin{equation}
\label{E:3.15q}
\|\nabla\wt f\|_{L^2(\BR^2)}=\|\nabla 
f\|_{\ell^2(\Z^2)}.
\end{equation}
\end{enumerate}
Here $\nabla f$ denotes the discrete gradient whereas $\nabla\tilde f$ denotes the continuous gradient.
\end{lemma}

\begin{proofsect}{Proof}
This is a simplified version of Lemma~3.3 from Biskup, Fukushima and K\"onig~\cite{BFK14} which itself is a version of Lemma~2.1 in Becker and K\"onig~\cite{Becker-Koenig}. See also van der Hofstad, K\"onig and M\"orters~\cite[Proposition~5.1]{HKM}.
\end{proofsect}

With this lemma in hand, the proof of the above proposition is quite straightforward:

\begin{proofsect}{Proof of Proposition~\ref{prop-4.5}}
Fix~$\cmss P\in\mathfrak P_L$ and let~$S$ be such that $\cmss P(S)=\cmss P$. Let $g\colon\Z^2\to\BR$ denote the principal eigenfunction of the (lattice) Laplacian in~$S$ with zero boundary conditions outside of~$S$. Thanks to the above definitions (and using~$g$ for function~$f$ in Lemma~\ref{fe}) there is a function $h\colon\BR^2\to\BR$ (corresponding to~$\widetilde f$ in the lemma) such that
\begin{enumerate}
\item[(1)] $h\in C^1_\cc(V_L)$,
\item[(2)] $\Vert\nabla h\Vert_{L^2(\BR^2)}=\Vert \nabla g\Vert_{\ell^2(\Z^2)}$,
\item[(3)] $\Vert h\Vert_{L^2(\BR^2)}\ge \Vert g\Vert_{\ell^2(\Z^2)}-C\Vert \nabla g\Vert_{\ell^2(\Z^2)}$.
\end{enumerate}
Hence, assuming that the right-hand side of the expression in (3) is positive,
\begin{equation}
\lambda(V_L)\le\frac{\Vert\nabla h\Vert_{L^2(\BR^2)}^2}{\Vert h\Vert_{L^2(\BR^2)}^2}
\le \frac{\Vert \nabla g\Vert_{\ell^2(\Z^2)}^2}{(\Vert g\Vert_{\ell^2(\Z^2)}-C\Vert \nabla g\Vert_{\ell^2(\Z^2)})^2}.
\end{equation}
Under the normalization $\Vert g\Vert_{\ell^2(\Z^2)}=1$ we have $\Vert \nabla g\Vert_{\ell^2(\Z^2)}=\sqrt{\lambda_S^{\ssup 1}}$ and the claim follows.
\end{proofsect}

The upshot of Proposition~\ref{prop-4.5} is that, once $\lambda^{\ssup1}_S$ is known to be small, we get a tight comparison between $\lambda_S^{\ssup1}$ and $\lambda(V_L)$. It remains relate $\CP(\Int(\cmss P))$ to $\CP(V_L)$. This is the content of:

\begin{proposition}
\label{prop-4.7}
For each norm~$\rho$ there is a constant $c\in(0,\infty)$ such that the following holds for the perimeter functional $U\mapsto\CP(U)$ defined using~$\rho$:
For each $L\in\BN$ and each $\cmss P\in\mathfrak P_L$, if $V_L$ is related to $\cmss P$ as in \eqref{E:UL}, then
\begin{equation}
\label{E:4.26r}
\CP\bigl(\Int(\cmss P)\bigr)\ge \CP(V_L)-c\,L.
\end{equation}
Moreover, for each $\epsilon>0$ there are $\delta>0$ and~$\zeta>0$ such that for any $r>0$ and any minimizing shape~$U_0\in\CU$ (i.e., a set with $\CF(U_0)=\min\CF$),
\begin{multline}
\label{E:4.29r}
\qquad
\distH\bigl(\hull(\cmss P),rU_0\bigr)<\delta r\quad\text{\rm and}\quad\distH\bigl(\hull(\cmss P),\Int(\cmss P)\bigr)>\epsilon r
\\\qquad\Rightarrow\qquad 
\CP\bigl(\Int(\cmss P)\bigr)\ge \CP(V_L)-c\,L+\zeta r
\qquad
\end{multline}
\end{proposition}

For the proof we will need to recall some geometric facts established by the present authors in a companion paper (Biskup and Procaccia~\cite{BP-analysis}). First, setting $U+V := \{x+y\colon x\in U, y\in V\}$ for any sets $U,V\subset\BR^2$, a particular feature of the two-dimensional perimeter functional is that
\begin{equation}
\label{E:3.39ue}
U,V\in\CU\quad\text{convex}\quad\Rightarrow\qquad \CP(U+V)=\CP(U)+\CP(V),
\end{equation}
see~\cite[Lemma]{BP-analysis}. This alone now permits us to give:

\begin{proofsect}{Proof of Proposition~\ref{prop-4.7}, formula \eqref{E:4.26r}}
Fix~$\cmss P$ and consider the perimeter with respect to the norm $\rho$. Note that the vertices of~$\hull(\cmss P)$ are also vertices of~$\cmss P$. The triangle inequality then readily shows
\begin{equation}
\CP\bigl(\hull(\cmss P)\bigr)\le\CP\bigl(\Int(\cmss P)\bigr).
\end{equation}
If $B:=\{x\in\BR^2\colon|x|_1<L+3\}$, then $V_L=\hull(\cmss P)+B$. By \eqref{E:3.39ue},
\begin{equation}
\label{E:3.44ue}
\CP(V_L)=\CP\bigl(\hull(\cmss P)\bigr)+\CP(B).
\end{equation}
Hereby \eqref{E:4.26r} follows by the fact that $\CP(B)\le cL$ holds for all~$L\ge1$ with some~$c\in(0,\infty)$ that only depends on the underlying norm~$\rho$.
\end{proofsect}

For the second part of Proposition~\ref{prop-4.7} we first note that, since the minimizer of~$\CF$ is unique up to translates (see~\cite[Theorem~1.1]{BP-analysis}), it suffices to consider just one minimizing shape~$U_0$. Here we recall that~$U_0$ is said to contain a facet in direction of a unit vector~$e\in\BR^2$, if there are $x,y$ with $y=x+|x-y|_2e$, such that the linear segment $[x,y]$ is contained entirely in~$\partial U_0$. If $e'$ is a unit vector that is orthogonal to~$e$ and define, in light of convexity $s\mapsto\rho(e+se')$,
\begin{equation}
\theta^\pm:=\frac{\textd}{\textd s^\pm}\rho(e+se')\Bigl|_{s=0}.
\end{equation}
The difference $\theta^+-\theta^-$ is non-negative and independent of the orientations of~$e$ and~$e'$ --- and thus depends only on~$e$. Then
\begin{equation}
\label{E:3.41ue}
U_0\text{ has a facet in direction }e\quad\Leftrightarrow\quad \theta^+>\theta^-,
\end{equation}
see \cite[Theorem~1.5]{BP-analysis}. The direction~$e$ with the property on the right is then called \emph{degenerate}.

The reason why we care for degenerate directions here is as follows: The $\rho$-length of a path that connects two points on the facet can then be compared to the $\rho$-length of the linear segment connecting just the endpoints. Indeed, for any pair of orthogonal unit vectors~$e$ and~$e'$, 
\begin{equation}
\label{E:2.40a}
\sum_{i=1}^n\rho(t_ie+u_i e')\ge\Bigl(\,\sum_{i=1}^n t_i\,\Bigr)\rho(e)+\frac{\theta^+-\theta^-}2\sum_{i=1}^n|u_i|
\end{equation}
holds for any $t_1,\dots,t_n\in\BR$ and any $u_1,\dots,u_n\in\BR$ such that $\sum_{i=1}^n u_i=0$; see \cite[Lemma~4.1]{BP-analysis}. To use these facts efficiently, we will also need:

\begin{lemma}
\label{lemma-3.7ue}
Let $U\in\CU$ be convex and let $B_r$ denote the Euclidean ball of radius~$r$. For each $\eta>0$ and each~$\xi>0$ there is $\delta>0$ such that if $x,y$ obey
\begin{equation}
[x,y]\subset\partial U+B_\delta,
\end{equation}
then either $|x-y|_2<\eta$ or there are $z,z'\in\partial U$ with
\begin{equation}
|z-x|_2<\xi,\quad |z'-y|_2<\xi\quad\text{and}\quad [z,z']\subset\partial U.
\end{equation}
\end{lemma}

\begin{proofsect}{Proof}
If the first alternative fails for each~$\delta>0$ existed then there would exist sequences $\{x_n\}$ and~$\{y_n\}$ with $x_n\to z\in\partial U$ and~$y_n\to z'\in\partial U$ such that $[x_n,y_n]\subset\partial U+B_{1/n}$ and yet $|x_n-y_n|_2\ge\eta$. By convexity of~$U$, the segment $[z,z']$ would then lie on a facet of~$U$.
\end{proofsect}

With these in hand, we can now give:

\begin{proofsect}{Proof of Proposition~\ref{prop-4.7}, formula \eqref{E:4.29r}}
Let~$c_1,c_2\in(0,\infty)$ be such that $c_1|x|_2\le\rho(x)\le c_2|x|_2$ for each~$x\in\BR^2$. Fix~$\epsilon>0$ and pick~$\eta>0$ such that $2c_1\epsilon>c_2\eta$. Next let~$m$ denote the minimum of $(\theta^+-\theta^-)/2$ for all degenerate directions~$e$ for which~$U_0$ has a facet of length at least~$\eta/2$. Since~$\partial U_0$ is rectifiable, there are only a finite number of such facets and so, by \eqref{E:3.41ue},~$m>0$.
Then let~$\xi>0$ be such that $m(\epsilon-\eta)>4c_2\xi$ and $\xi<\eta/4$. For these~$\eta$ and~$\xi$, let now~$\delta>0$ be such that Lemma~\ref{lemma-3.7ue} applies. We claim that if~$\cmss P$ be such that the two conditions on the left of \eqref{E:4.29r} hold for some~$r>0$ and the~$\epsilon$ and~$\delta$ as above, then
 \begin{equation}
\label{E:3.47ue}
\CP(\Int(\cmss P))\ge\CP(\hull(\cmss P))+\zeta r
\end{equation}
with
\begin{equation}
\label{E:3.48ue}
\zeta:=\min\bigl\{(2c_1\epsilon-c_2\eta),m(\epsilon-\eta)-4c_2\xi\bigr\}
\end{equation}
which is positive by our choices above. This is enough to \eqref{E:4.29r} by invoking the argument after \eqref{E:3.44ue}. It thus remains to show that \eqref{E:3.47ue} indeed holds.

We begin by noting that, since $\Int(\cmss P)\subseteq\hull(\cmss P)$, the assumption $\distH(\hull(\cmss P),\Int(\cmss P))>\epsilon r$ implies the existence of a vertex~$v\in\cmss P$ such that
\begin{equation}
\dist_2\bigl(v,\hull(\cmss P)^\cc\bigr)\ge\epsilon r.
\end{equation}
As $\hull(\cmss P)$ is a polygonal domain whose every vertex is a vertex of~$\cmss P$, there are vertices $x,y$ of~$\hull(\cmss P)$ such that the part of~$\cmss P$ between~$x$ and~$y$ (in the chosen orientation) passes through~$v$. These can in fact be chosen so that $[x,y]\subset\partial\hull(\cmss P)$. By the triangle inequality,
\begin{equation}
\CP(\Int(\cmss P))\ge\CP(\hull(\cmss P))+\rho(v-x)+\rho(v-y)-\rho(x-y).
\end{equation}
We will now bound the expression involving norms on the right-hand side. 

Let $\BH_{xy}$ denote the open half-plane containing~$v$ whose boundary line passes through~$x$ and~$y$. Then $\hull(\cmss P)\subset\BH_{xy}$ and so
\begin{equation}
\dist_2(v,\BH_{xy}^\cc)\ge \dist_2\bigl(v,\hull(\cmss P)^\cc\bigr)\ge\epsilon r.
\end{equation}
Thus, in particular,
\begin{equation}
|v-x|_2\ge\epsilon r\quad\text{and}\quad |v-y|_2\ge\epsilon r.
\end{equation}
Next, by our assumptions,
\begin{equation}
\dist_2\bigl(x,\partial (rU_0)\bigr)<\delta r\quad\text{and}\quad\dist_2\bigl(y,\partial (rU_0)\bigr)<\delta r
\end{equation}
and so $[x,y]\subset (rU_0)+B_{r\delta}$.
Our choice of~$\delta$ (and scale invariance) then ensures that then either $|x-y|_2<\eta r$ or there is a segment $[z,z']$ on $\partial(rU_0)$ such that $|x-z|_2<\xi r$ and $|y-z'|_2<\xi r$. In the former case we have
\begin{equation}
\rho(x-y)\le c_2\eta r\quad\text{while}\quad\rho(v-x)+\rho(v-y)\ge 2c_1\epsilon r
\end{equation}
and so \eqref{E:3.47ue} holds with the first alternative in \eqref{E:3.48ue}.
In the latter case we add the segment $[z,x]$ and~$[y,z']$ to the path from~$x$ to~$z$ to~$y$ and apply \eqref{E:2.40a} with the result
\begin{equation}
\rho(x-v)+\rho(y-v)+\rho(x-z)+\rho(y-z')\ge\rho(z-z')+m(\epsilon-\eta)r.
\end{equation}
Bounding $\rho(x-z)+\rho(y-z')\le 2c_2r\xi$ and $\rho(z-z')\ge\rho(x-y)-2c_2r\xi$, we thus get \eqref{E:3.47ue} with the second alternative in \eqref{E:3.48ue}. Hence \eqref{E:3.47ue} is proved and the claim holds.
\end{proofsect}

\section{Proof of Main Theorems}
\label{sec5}\nopagebreak\noindent
We are now ready to move to the proof of the shape theorem. Throughout this section $\rho$ denotes the norm constructed in Theorem~\ref{thm-norm} and~$\CF$ is the functional in \eqref{E:1.13} with the perimeter defined using this norm. The starting point is a lower bound on the partition function $Z(t,\beta)$. This will set a scale to which we will later compare the contribution of various undesirable events.

\subsection{Large deviation lower bound}
The sole subject of this subsection is the proof of:

\begin{proposition}
\label{prop-5.1}
We have
\begin{equation}
\liminf_{\beta\to\infty}\,\,\liminf_{t\to\infty}\,\,\frac{r(t,\beta)^2}t\,\log\,Z(t,\beta)
\ge-\min\CF.
\end{equation}
\end{proposition}

For the proof we will need:

\begin{lemma}
\label{lemma-5.2}
Given $U\in\CU$, let~$S=S(U,t,\beta)$ be as in \eqref{E:3.12}. For each $\epsilon>0$ there is~$r_0=r_0(U)$ such that if $r(t,\beta)\ge r_0$, then
\begin{equation}
\lambda_S^{\ssup1}\le(1+\epsilon)\frac{\lambda(U)}{r(t,\beta)^2}.
\end{equation}
\end{lemma}

\begin{proofsect}{Proof}
Let $U\in\CU$, let $S=S(U,t,\beta)$ and abbreviate $r:=r(t,\beta)$. Fix~$\epsilon>0$ and pick $h\in C_c^\infty(U)$ such that $\Vert h\Vert_{L^2(\BR^2)}=1$ and $\Vert\nabla h\Vert_{L^2(\BR^2)}^2\le(1+\epsilon)^{1/2}\lambda(U)$. Then define
\begin{equation}
\label{E:5.3a}
g(x):=\int_{[0,1]^2}h\bigl(r^{-1}(x+z)\bigr)\,\,\textd z,\qquad x\in\BZ^2
\end{equation}
and note that, for~$r$ sufficiently large, we have $\supp(g)\subset S$. We claim that
\begin{equation}
\label{E:5.4}
\Vert\nabla g\Vert_{\ell^2(\Z^2)}\le\frac1{r}\Vert\nabla h\Vert_{L^2(\BR^2)}
\end{equation}
and, for some absolute constant~$c\in(0,\infty)$, also
\begin{equation}
\label{E:5.5}
0\le\Vert h\Vert_{L^2(\BR^2)}^2-\Vert g\Vert_{\ell^2(\Z^2)}^2\le \frac c{r^2}\Vert\nabla h\Vert_{L^2(\BR^2)}^2.
\end{equation}
These then readily give the claim via
\begin{equation}
\lambda_S^{\ssup 1}\le\frac{\Vert\nabla g\Vert_{\ell^2(\Z^2)}^2}{\Vert g\Vert_{\ell^2(\Z^2)}^2}
\le\frac1{r^2}\frac{\Vert\nabla h\Vert_{L^2(\BR^2)}^2}{\Vert h\Vert_{L^2(\BR^2)}^2-\frac c{r^2}\Vert\nabla h\Vert_{L^2(\BR^2)}^2}\le\frac1{r^2}\,\frac{(1+\epsilon)^{1/2}\lambda(U)}{1-\frac{c'}{r^2}\lambda(U)}
\end{equation}
where $c':=c(1+\epsilon)^{1/2}$ and where we assumed that $\frac{c'}{r^2}\lambda(U)<1$. Indeed, increasing~$r$ further if necessary, the right-hand side is at most $r^{-2}(1+\epsilon)\lambda(U)$.

It thus remains to establish the bounds \twoeqref{E:5.4}{E:5.5}. The first of these is quite straightforward: Letting $e_1$ and~$e_2$ denote the unit vectors in coordinate directions, we have
\begin{equation}
g(x+e_i)-g(x) = \frac1{r}\int_{[0,1]^2}\textd z\int_0^1\textd s\,\,e_i\cdot\nabla h\bigl(r^{-1}(x+z+se_i)\bigr).\end{equation}
Applying Cauchy-Schwarz, we then get
\begin{equation}
\Vert\nabla g\Vert_{\ell^2(\Z^2)}^2
\le\frac1{r^2}\sum_{x\in\Z^2}\sum_{i=1,2}\int_{[0,1]^2}\textd z\int_0^1\textd s\,\Bigl|e_i\cdot\nabla h\bigl(r^{-1}(x+z+se_i)\bigr)\Bigr|^2.
\end{equation}
As is now easy to check, the sums on the right-hand side then reduce to $\Vert\nabla h\Vert_{L^2(\BR^2)}^2$.

The second inequality \eqref{E:5.5} is slightly more involved. First, the bound on the left is obtained by using Cauchy-Schwarz in \eqref{E:5.3a} and summing over~$x\in\Z^2$. This in turn shows
\begin{equation}
\begin{aligned}
\Vert h\Vert_{L^2(\BR^2)}^2-&\Vert g\Vert_{\ell^2(\Z^2)}^2
\\&=\sum_{x\in\Z^2}\biggl(\int_{[0,1]^2}h\bigl(r^{-1}(x+z)\bigr)^2\,\textd z-\Bigl(\int_{[0,1]^2}h\bigl(r^{-1}(x+z)\bigr)\,\textd z\Bigr)^2\biggr).
\end{aligned}
\end{equation}
A simple rewrite and the Cauchy-Schwarz estimate then permit us to recast this as
\begin{equation}
\begin{aligned}
\sum_{x\in\Z^2}&\int_{[0,1]^2\times[0,1]^2}\textd z\,\textd z'
\,\,\Bigl|\,h\bigl(r^{-1}(x+z)\bigr)-h\bigl(r^{-1}(x+z')\Bigr|^2\\
&\le\frac1{r^2}\sum_{x\in\Z^2}\int_0^1\textd s \int_{[0,1]^2\times[0,1]^2}\textd z\,\textd z'\,\,\,
\biggl|(z-z')\cdot\nabla h\Bigl(r^{-1}\bigl(x+sz+(1-s)z'\bigr)\Bigr)\biggr|^2.
\end{aligned}
\end{equation}
The expression in absolute value is bounded by twice the norm of~$\nabla h$ at the corresponding point.
The push-forward on~$\BR^2$ of the integrating measure under the map $(s,z,z')\mapsto sz+(1-s)z'$ is dominated by a constant times the Lebesgue measure on, say, $[-1,2]\times[-1,2]$. Using translation invariance again, we get the inequality on the right of \eqref{E:5.5} as well.
\end{proofsect}

We will also need:

\begin{lemma}
\label{lemma-4.3ue}
For each admissible~$\BP$ there is~$q\ge1$ such that for $\BP$-a.e.\ sample of edge weights, each~$x,y\in\BZ^2$ and each path $\gamma=(x_0^\star,\dots,x_n^\star)\in\Gamma(x,y)$ with $\frakd(\gamma)\le 2D(x,y)$ we have
\begin{equation}
\max_{i=0,\dots,n}\,\max\bigl\{|x_i^\star-x|_2,|x_i^\star-y|_2\bigr\}\le q|x-y|_2.
\end{equation}
\end{lemma}

\begin{proofsect}{Proof}
Recall that for each admissible~$\BP$ there are $0<\kappa<\upsilon<\infty$ such that the edge weights lie in $[\kappa,\upsilon]$ $\BP$-a.s. Let $c,c'\in(0,\infty)$ be the constant such that $c|x|_2\ge|x|_1\ge c'|x|_2$ for all~$x\in\BR^2$. Clearly,
\begin{equation}
\frakd(\gamma)\ge a\max_{i=0,\dots,n}\,\max\bigl\{|x_i^\star-x_0^\star|_1,|x_i^\star-x_n^\star|_1\bigr\}
\end{equation}
so if the maximum in the statement is larger than $r|x-y|_2$ then we have $\frakd(\gamma)\ge \kappa(c'/c)q|x-y|_1$. But the assumptions also tell us $\frakd(\gamma)\le2D(x,y)\le 2\upsilon |x-y|_1$ and so this is not possible once $q$ is so large that $2\upsilon<\kappa (c'/c)q$.
\end{proofsect}

\begin{proofsect}{Proof of Proposition~\ref{prop-5.1}}
Let~$U_0$ be the minimizer of~$\CF$ and let $\epsilon>0$ be so small that Proposition~\ref{prop-3.2} applies to $U:=(1+\epsilon)U_0$. Let~$\cmss P=(x_0,\dots,x_n)$ be a polygonal curve such that $0\in\Int(\cmss P)\subset U_0$, the vertices of~$\cmss P$ lie on~$\partial U_0$ and, for~$\delta:=\min\{|x_i-x_{i-1}|_2\colon i=1,\dots,n\}$, we have
\begin{equation}
\label{E:4.13ue}
(1-\epsilon)U_0\cap\bigl(\cmss P+B_{2q\delta}\bigr)=\emptyset\quad\text{and}\quad
\bigl(\cmss P+B_{2q\delta}\bigr)\subseteq(1+\epsilon)U_0,
\end{equation}
where~$B_R$ is the Euclidean ball of radius~$R$ centered at the origin and~$q$ is as in Lemma~\ref{lemma-4.3ue}.
This is possible in light of convexity of~$U_0$. 
 
Next abbreviate~$r:=r(t,\beta)$ and assume that~$r$ is so large that 
\begin{equation}
\label{E:4.14ue}
D\bigl(\lfloor rx_i\rfloor,\lfloor rx_{i-1}\rfloor\bigr)\le (1+\epsilon)r\rho(x_i-x_{i-1}),\qquad i=1,\dots,n.
\end{equation}
This is possible thanks to Theorem~\ref{thm-norm}, the fact that~$\cmss P$ is bounded and also thanks to the continuity of~$x\mapsto\rho(x)$. Pick a minimizing path $\gamma_i\in\Gamma(\lfloor rx_i\rfloor,\lfloor rx_{i-1})$ for each $i=1,\dots,n$. These paths form a closed cycle on~$\BZ^2$ which, however, may not be simple. Notwithstanding, letting~$S$ denote the finite connected component of the interior of the above cycle containing the origin, Lemma~\ref{lemma-4.3ue} and \eqref{E:4.13ue} ensure that
\begin{equation}
\tilde S\subseteq S\subseteq S(\bigl(1-\epsilon)U_0,t,\beta)
\quad\text{for}\quad \tilde S:=S\bigl((1-\epsilon)U_0,t,\beta)
\end{equation}
and thus, in particular, $S\in\mathfrak S$.  

The inclusion $\tilde S\subseteq S$ and Lemma~\ref{lemma-5.2} now yield
\begin{equation}
\lambda^{\ssup1}_S\le \lambda^{\ssup1}_{\tilde S}\le \frac{1+\epsilon}{(1-\epsilon)^2}\frac{\lambda(U_0)}{r(t,\beta)^2}.
\end{equation}
On the other hand, the construction of~$S$ and \eqref{E:4.14ue} guarantee
\begin{equation}
\begin{aligned}
H(S)\le\sum_{i=1}^n\frakd(\gamma_i)
&=\sum_{i=1}^n D\bigl(\lfloor rx_i\rfloor,\lfloor rx_{i-1}\rfloor\bigr)\\
&\le(1+\epsilon)r\sum_{i=1}^n\rho(x_i-x_{i-1})\le(1+\epsilon)r\CP(U_0).
\end{aligned}
\end{equation}
Since $|\partial S|\le a^{-1}H(S)\le c r\CP(U_0)$ where $a$ is the lower bound on the weights in~$\BP$, Proposition~\ref{prop-3.2} implies
\begin{equation}
Z(t,\beta)\ge  \texte^{-t\lambda^{\ssup1}_S-\beta H(S)-c|\partial S|}
\end{equation}
as soon as $\beta$ and~$t$ are sufficiently large.
From \eqref{E:1.9ue} we thus get
\begin{equation}
\begin{aligned}
\log Z(t,\beta)
&\ge-\frac t{r^2}\frac{1+\epsilon}{(1-\epsilon)^2}\lambda(U_0)-(\beta+ca^{-1})(1+\epsilon)r\CP(U_0)
\\
&\ge-\frac t{r^2}\frac{1+\epsilon}{(1-\epsilon)^2}\Bigl(\CF(U_0)- \frac c{a\beta}\CP(U_0)\Bigr).
\end{aligned}
\end{equation}
Taking the limits $t\to\infty$, $\beta\to\infty$ and $\epsilon\downarrow0$, the claim follows.
\end{proofsect}

\subsection{Proof of the shape theorem}
Suppose, throughout this subsection, that~$\kappa \in(0,1)$ is an $\BP$-a.s.\ lower bound on edge weights.
Before we delve into proof of the shape theorem, we begin with some basic lemmas. The first one notes a restriction on the diameter of~$R(t)$:

\begin{lemma}
\label{lemma-diameter}
For each $\beta>\frac2 \kappa \log 4$, all $M>0$ and all~$t$ sufficiently large,
\begin{equation}
E^0\Bigl(\, \texte^{-\beta H(R(t))}\1_{\{|\partial R(t)|>M r(t,\beta)\}}\Bigr)
\le  \texte^{-\frac12 \kappa M\, t/r(t,\beta)^2}.
\end{equation}
\end{lemma}

\begin{proofsect}{Proof}
Note that $H(R(t))\ge \kappa|\partial R(t)|$. On the stated event, $\partial R(t)$ can be identified with a closed path on the dual-$\Z^2$ of length at least $M r(t,\beta)$ surrounding the origin in~$\Z^2$. The number of such paths of length~$n$ is at most $n4^n$ and so for $\beta_0:=\frac1\kappa\log4$ and $\beta>\beta_0$,
\begin{equation}
\label{E:5.3}
E^0\Bigl(\, \texte^{-\beta H(R(t))}\1_{\{|\partial R(t)|>M r(t,\beta)\}}\Bigr)
\le\sum_{n>M r(t,\beta)}n(4 \texte^{-\beta \kappa})^{n}\le C(\beta) \texte^{-(\beta-\beta_0)\kappa M \,r(t,\beta)}.
\end{equation}
Now for $\beta>2\beta_0$ we have $(\beta-\beta_0)>\frac12\beta$ and, by definition of~$r(t,\beta)$, 
\begin{equation}
(\beta-\beta_0) r(t,\beta)>\frac12 \frac t{r(t,\beta)^2}.
\end{equation}
For~$t$ sufficiently large, we can then absorb~$C(\beta)$ into the exponent. 
\end{proofsect}

The next lemma effectively restricts $\lambda_{R(t)}^{\ssup1}$ to values order~$r(t,\beta)^{-2}$:

\begin{lemma}
\label{lemma-eigenvalue}
For each $\beta>\frac 3\kappa\log 4$ and all~$M$ and~$t$ sufficiently large,
\begin{equation}
E^0\Bigl(\, \texte^{-\beta H(R(t))}\1_{\{\lambda_{R(t)}^{\ssup1}>M r(t,\beta)^{-2}\}}\Bigr)
\le 2 \texte^{-\frac12\kappa M \,t/r(t,\beta)^2}.
\end{equation}
\end{lemma}

\begin{proofsect}{Proof}
By Lemmas~\ref{lemma-diameter} and~\ref{prop-3.1w},
\begin{multline}
\qquad
E^0\Bigl(\, \texte^{-\beta H(R(t))}\1_{\{\lambda_{R(t)}^{\ssup1}>M r(t,\beta)^{-2}\}}\Bigr)
\le \texte^{-\frac12\kappa M\,t/r(t,\beta)^2}
\\
+\sum_{S\in\mathfrak S}
\1_{\{|\partial S|\le Mr(t,\beta)\}}\1_{\{\lambda_S^{(1)}>Mr(t,\beta)^{-2}\}}
|S|^{3/2} \texte^{-\beta H(S)-t\lambda_S^{\ssup1}}
\qquad
\end{multline}
Bounding $|S|^{3/2}$ by $c M^{3/2}r(t,\beta)^3$ for some constant $c\in(0,\infty)$ and employing the restriction on~$\lambda_S^{\ssup1}$ dominates the second term on the right by
\begin{equation}
cM^{3/2}r(t,\beta)^3 \texte^{-Mt/r(t,\beta)^2}\,\sum_{S\in\mathfrak S} \texte^{-\beta H(S)}.
\end{equation}
For $\beta>\beta_0:=\frac1\kappa\log 4$, the sum is bounded by a constant that only depends on~$\beta$. Taking~$M$ sufficiently large, the result is thus at most another factor of~$ \texte^{-\frac12\kappa M\,t/r(t,\beta)^2}$ (since $\kappa<1$).
\end{proofsect}

We will now use these to prove:

\begin{proposition}
\label{prop-5.5}
Let $\CM:=\{x+U_0\colon -x\in U_0\}$ where $U_0$ is a minimizer of~$\CF$. For each~$\epsilon>0$ there is~$\delta>0$ such that
\begin{equation}
\limsup_{\beta\to\infty}\,\,\limsup_{t\to\infty}\frac{r(t,\beta)^2}t\,\log E^0\Bigl(\, \texte^{-\beta H(R(t))}\1_{\{\distH(r(t,\beta)^{-1}R(t),\CM)>\epsilon\}}\Bigr)\le-\min\CF-\delta.
\end{equation}
\end{proposition}

\begin{proofsect}{Proof}
Much of the work has already been done; here we just basically assemble all facts together. Abbreviate $r:=r(t,\beta)$ and let $L:=\lfloor 2\zeta r\rfloor$ for some~$\zeta>0$ to be let go to zero at the end. Let
\begin{equation}
\mathfrak S_{M,\zeta}:=\big\{S\in\mathfrak S\colon 2L<|\partial S|\le M r,\,\lambda_{S}^{\ssup1}\le M r^{-2}\bigr\}.
\end{equation}
Since each $S\in\mathfrak S_{M,\zeta}$ is connected, the discrete Faber-Krahn estimate $\lambda_S^{\ssup1}\ge c|S|^{-1}$ along with the isoperimetric inequality $|\partial S|\ge c'|S|^{1/2}$ (both written for $d=2$) show that~$|\partial S|>2L$ is implied by~$\lambda_{S}^{\ssup1}\le Mr^{-2}$ as soon as $4\zeta\sqrt M<c'\sqrt c$.
Thanks to Lemmas~\ref{lemma-diameter}--\ref{lemma-eigenvalue}, we can choose~$M$ so that for all~$\zeta$ sufficiently small the event $\{R(t)\in\mathfrak S_{M,\zeta}\}$ may freely be inserted into the expectation. This reduces the proof to a suitable estimate on the quantity
\begin{equation}
Y_\epsilon:=\sum_{\begin{subarray}{c}
S\in\mathfrak S_{M,\zeta}\\\distH(r^{-1}S,\CM)>4\epsilon
\end{subarray}}
 \texte^{-\beta H(S)}P^0\bigl(R(t)=S\bigr).
\end{equation}
where we write~$4\epsilon$ instead of~$\epsilon$ for later convenience.

We start by applying the inequality from Lemma~\ref{prop-3.1w}, use that $|S|\le c Mr^2$ for some constant $c\in(0,\infty)$ and then represent (as an upper bound) the sum over~$S$ as a sum over $L$-skeletons and the sum over~$S$ for a given skeleton. This yields
\begin{equation}
\begin{aligned}
Y_\epsilon
&\le
\sum_{\begin{subarray}{c}
S\in\mathfrak S_{M,\zeta}\\\distH(r^{-1}S,\CM)>4\epsilon
\end{subarray}}
 \texte^{-\beta H(S)}|S|^{3/2} \texte^{-t\lambda_S^{\ssup1}}
\\
&\le \,\,\,\,(cM)^{3/2}r^3\!\!\!\!\!\sum_{\begin{subarray}{c}
\cmss P\in\mathfrak P_L\\\distH(r^{-1}\cmss P,\CM)>2\epsilon
\end{subarray}}\,\,
\,\sum_{\begin{subarray}{c}
S\in\mathfrak S_{M,\zeta}\\\cmss P(S)=\cmss P
\end{subarray}}
 \texte^{-\beta H(S)} \texte^{-t\lambda_S^{\ssup1}},
\end{aligned}
\end{equation}
where we invoked the bound $\distH(S,\Int(\cmss P))\le L\le 2\epsilon r$ (which requires $\zeta<\epsilon$) along with $\distH(\Int(\cmss P),\CM)=\distH(\cmss P,\CM)$ to move the Hausdorff-distance restriction to~$\cmss P$. We also used that $|\partial S|>2L$ implies $\cmss P(S)\in\mathfrak P_L$.

Let now $V_L$ be related to~$\cmss P$ as in~\eqref{E:UL}. Proposition~\ref{prop-4.5} then shows
\begin{equation}
\lambda_S^{\ssup1}\ge\lambda(V_L)\bigl(1-CM^{1/2}r^{-1}\bigr)^2\ge\lambda(V_L)(1-\zeta)
\end{equation}
for every $S\in\mathfrak S_{M,\zeta}$ once~$r$ is large enough. Denote
\begin{equation}
\mathfrak P_L':=\bigl\{\cmss P(S)\colon S\in\mathfrak S_{M,\zeta}\bigr\},
\end{equation}
and note that every $\cmss P\in\mathfrak P_L'$ has at most $Mr/L+1$ linear segments. This is less than $1/\zeta$ for~$r$ large enough so bounding the sum over~$S$ with the help of Lemma~\ref{lemma-4.4} yields
\begin{equation}
\label{E:5.25k}
Y_\epsilon\le 
\,\,\,\,(cM)^{3/2}r^3 \texte^{2b\beta L}\!\!\!\!\!\sum_{\begin{subarray}{c}
\cmss P\in\mathfrak P_L'\\\distH(r^{-1}\cmss P,\CM)>2\epsilon
\end{subarray}}
 \texte^{-t\lambda(V_L)(1-\zeta)-\beta(1-\zeta)\CP(\Int(\cmss P))}.
\end{equation}
Our next goal is to invoke Proposition~\ref{prop-4.7} so let~$\delta$ correspond to~$\epsilon$ in the statement. Without loss of generality we may assume that~$\delta<\epsilon$. Then for each~$\cmss P$ subject to $\distH(r^{-1}\cmss P,\CM)>2\epsilon$,
\begin{enumerate}
\item[(1)] either $\distH(r^{-1}\hull(\cmss P),\CM)\ge\delta$
\item[(2)] or $\distH(r^{-1}\hull(V_L),\CM)<\delta$ and $\distH\bigl(\Int(\cmss P),\hull(\cmss P)\bigr)>\epsilon\,r$.
\end{enumerate}
In the former case we invoke \eqref{E:4.26r} while in the latter case we use \eqref{E:4.29r} to conclude, for some numerical constants~$c_1,c_2\in(0,\infty)$ that
\begin{multline}
\qquad
\text{sum in \eqref{E:5.25k}}
\le \texte^{c_1\beta L}\sum_{\begin{subarray}{c}
\cmss P\in\mathfrak P_L'\\\distH(r^{-1}\hull(\cmss P),\CM)>\delta
\end{subarray}}
 \texte^{-(1-\zeta)\frac t{r^2}\CF(V_L(\cmss P))}
\\
+ \texte^{c_1\beta L-c_2\beta\epsilon r}\sum_{\cmss P\in\mathfrak P_L'}
 \texte^{-(1-\zeta)\frac t{r^2}\CF(V_L(\cmss P))}
\qquad
\end{multline}
Since $\distH(r^{-1}\hull(\cmss P),\CM) = \distH(r^{-1}V_L,\CM)-L/r$, once~$\zeta$ is so small that $\zeta<\delta/2$, \eqref{E:1.15} then implies the existence of $\beta(\epsilon)<\infty$ and~$\delta'>0$ such that for all $\beta>\beta(\epsilon)$,
\begin{equation}
\distH(r^{-1}\hull(\cmss P),\CM)>\delta
\quad\Rightarrow\quad
\CF\bigl(V_L(\cmss P)\bigr)\ge\min\CF+\delta'.
\end{equation}
Since $\beta r = \frac t{r^2}$, assuming without loss of generality that $\delta'\le c_2\epsilon$ yields, for some $c_3\in(0,\infty)$,
\begin{equation}
Y_\epsilon\le 2(cM)^{3/2}r^3\, \texte^{\beta c_3 L}\, \texte^{-\frac t{r^2}(1-\zeta)[\min\CF+\delta']}\,|\mathfrak P_L'|.
\end{equation}
Since every $\cmss P\in\mathfrak P_L'$ has at most $1/\zeta$ linear segments each of which is of length at most~$L$, there is~$c\in(1,\infty)$ such that
\begin{equation}
|\mathfrak P_L'|\le \bigl(c L^2\bigr)^{1/\zeta}.
\end{equation}
The entropy of the skeletons is thus negligible and so, since $\beta L=(1+o(1))2\zeta t/r^2$, we get
\begin{equation}
\beta>\beta(\epsilon)\quad\Rightarrow\quad
\limsup_{t\to\infty}\,\frac{r(t,\beta)^2}t\log Y_\epsilon\le -(1-\zeta)[\min\CF-\delta']-2c_3\zeta.
\end{equation}
Taking $\beta\to\infty$ followed by~$\zeta\downarrow0$ then gives the claim for all~$\epsilon>0$.
\end{proofsect}

It now remains to put the pieces together to get:

\begin{proofsect}{Proof of Theorem~\ref{thm-1}}
By Propositions~\ref{prop-5.1} and \ref{prop-5.5}, given~$\epsilon>0$, the probability in the statement decays exponentially in~$t/r(t,\beta)^{2}$  as soon as~$\beta$ is chosen sufficiently large. This readily yields the claim.
\end{proofsect}

\section{Laplacian eigenvalues and eigenfunctions}
\label{sec5a}\nopagebreak\noindent
The last item to finish is the proof of certain estimates for the eigenvalues and eigenfunctions of the Laplacian in discrete sets that well approximate a continuum domain. These have been deferred to here from Section~\ref{sec3}. The estimates are fairly standard; we include proofs for completeness of exposition. We begin with the claims concerning the eigenvalues:

\begin{proofsect}{Proof of Lemma~\ref{lemma-2.3}, \twoeqref{E:2.13ue}{E:3.18a}}
The inequality \eqref{E:2.13ue} follows either by Lemma~\ref{lemma-5.2} (and the corresponding statement in the continuum) or directly by the monotonicity of $S\mapsto\lambda^{\ssup1}_S$ and the fact that $S$ contains a box of side of order~$r(t,\beta)$. 

To get the spectral gap estimate \eqref{E:3.18a}, we have to work a bit harder. Fix~$U\in\CU$ and let~$\lambda^{\ssup k}(U)$ denote the $k$-th lowest eigenvalue of the negative Laplacian in~$U$ with Dirichlet boundary conditions on~$\partial U$. The argument will be based on the fact that
\begin{equation}
\Lambda^{\ssup2}(U):= \lambda^{\ssup1}(U)+\lambda^{\ssup2}(U)
\end{equation}
 obeys
\begin{equation}
\label{E:5.2a}
\Lambda^{\ssup2}(U)=\inf\bigl\{\Vert\nabla g_1\Vert_{L^2(\BR^2)}^2+\Vert\nabla g_2\Vert_{L^2(\BR^2)}^2\colon g_1,g_2\in C^\infty_\cc(U),\,\langle g_i,g_j\rangle_{L^2(\BR^2)}=\delta_{ij}\bigr\}.
\end{equation}
(This is sometimes called the Ky Fan principle; cf Ky Fan~\cite{KyFan}.)  For the discrete problem and
\begin{equation}
\Lambda_S^{\ssup2}:=\lambda^{\ssup1}_S+\lambda^{\ssup2}_S
\end{equation}
we similarly have
\begin{equation}
\label{E:5.4a}
\Lambda^{\ssup2}_S=\inf\bigl\{\Vert\nabla g_1\Vert_{\ell^2(\BZ^2)}^2+\Vert\nabla g_2\Vert_{\ell^2(\BZ^2)}^2\colon \supp(g_i)\subseteq S,\,\langle g_i,g_j\rangle_{\ell^2(\BZ^2)}=\delta_{ij}\bigr\}.
\end{equation}
Our argument is now based on the fact that the principal eigenvalue in any~$U\in\CU$ is non-degenerate, i.e.,  $\Lambda^{\ssup2}(U)>\lambda(U)$. Lemma~\ref{lemma-5.2} gives us a tight upper bound of~$\lambda^{\ssup1}_S$ in terms of~$\lambda^{\ssup1}(U)$; it thus suffices to show that, for each~$\epsilon>0$, 
\begin{equation}
\label{E:5.5a}
\Lambda^{\ssup2}_S\ge(1-\epsilon)\frac{\Lambda^{\ssup2}(U)}{r(t,\beta)^2}
\end{equation}
as soon as~$t$ and~$\beta$ are sufficiently large and~$S$ is any set as specified in the statement. As $S\mapsto\Lambda^{\ssup2}_S$ is non-decreasing with respect to inclusion and the scaling relation
\begin{equation}
\label{E:scale}
\lambda^{\ssup k}(\alpha U)=\alpha^{-2}\lambda^{\ssup k}(U)
\end{equation}
holds, it will in fact suffice to show this $S:=S(U,t,\beta)$ for any~$\epsilon>0$.  

The infimum in \eqref{E:5.4a} is achieved by $g_1:=h^{\ssup1}$ and $g_2:=h^{\ssup2}$, the first two eigenfunctions of the discrete Laplacian in~$S$. Recall that these are normalized to have $\ell^2(\Z^2)$-norm one. Let~$\widetilde g_1$, resp.,~$\widetilde g_2$ be the functions in~$\BR^2$ that are the counterparts to~$g_1$, resp.,~$g_2$ as in Lemma~\ref{fe}. These functions are not necessarily normalized or orthogonal. However, thanks to the fact that $\Vert \nabla g_i\Vert_{\ell^2(\BZ^2)}^2\le\Lambda^{\ssup2}_S=O(r(t,\beta)^{-2})$, \eqref{E:3.16q} shows, for some $c_1\in(0,\infty)$,
\begin{equation}
\Bigl|\Vert \widetilde g_i\Vert_{L^2(\BR^2)}-1\Bigr|\le c_1 r(t,\beta)^{-2},\qquad i=1,2.
\end{equation}
Since the map $f\mapsto\widetilde f$ in Lemma~\ref{fe} is linear, a similar argument applied to $f_\pm:=g_1\pm g_2$ combined with the polarization identity reveals that, for some $c_2\in(0,\infty)$,
\begin{equation}
\bigl|\langle \widetilde g_1,\widetilde g_2\rangle_{L^2(\BR^2)}\bigr|\le c_2 r(t,\beta)^{-2}
\end{equation}
Defining, for~$r(t,\beta)$ sufficiently large,
\begin{equation}
h_1:=\frac{\widetilde g_1}{\Vert \widetilde g_1\Vert_{L^2(\BR^2)}}
\quad\text{and}\quad
h_2:= \frac{\widetilde g_1-\langle h_1,\widetilde g_2\rangle_{L^2(\BR^2)}h_1}{\Vert \widetilde g_1-\langle h_1,\widetilde g_2\rangle_{L^2(\BR^2)}h_1\Vert_{L^2(\BR^2)}}
\end{equation}
we now get a pair of orthonormal functions. By \eqref{E:3.15q} and the fact that $g_1,g_2$ are normalized in~$\ell^2(\BZ^2)$ we then have, for some $c_3,c_4\in(0,\infty)$,
\begin{equation}
\Vert \nabla h_1\Vert_{L^2(\BR^2)}\le\bigl(1+c_3 r(t,\beta)^{-2}\bigr)\Vert \nabla g_2\Vert_{\ell^2(\BZ^2)}
\end{equation}
and
\begin{equation}
\Vert \nabla h_2\Vert_{L^2(\BR^2)}\le\bigl(1+c_3 r(t,\beta)^{-2}\bigr)\Vert \nabla g_2\Vert_{\ell^2(\BZ^2)}+c_4 r(t,\beta)^{-2}\Vert\nabla g_1\Vert_{\ell^2(\BZ^2)}.
\end{equation}
Since by \eqref{E:3.13aa} the support of~$h_1,h_2$ is contained in~$(1+\epsilon)U$, from \eqref{E:5.2a} and \eqref{E:5.4a} we get
\begin{equation}
\Lambda^{\ssup2}\bigl((1+\epsilon)r(t,\beta)U\bigr)\le \bigl(1+c_5 r(t,\beta)^{-2}\bigr)\Lambda^{\ssup2}_S
\end{equation}
for some $c_5\in(0,\infty)$.
Invoking the scaling relation \eqref{E:scale}, the bound \eqref{E:5.5a} follows.
\end{proofsect}

Next we move to the claims dealing with lower bounds on the principal eigenfunction:

\begin{proofsect}{Proof of Lemma~\ref{lemma-2.3}, \twoeqref{E:3.19}{E:3.17b}} The proof of \eqref{E:3.19} will be based on the following simple fact (derived, e.g., in Biskup and K\"onig~\cite[Lemma~4.1]{Biskup-Koenig}): Let~$Y_1,Y_2,\dots$ be the simple symmetric random walk on~$\BZ^2$ and, abusing our earlier notation slightly, let $\tau_S$ be the first exit time of the walk from~$S$. For the eigenfunction $h^{\ssup1}$ in~$S$, the process $\{M_{n\wedge\tau_S}\colon n\ge0\}$ where
\begin{equation}
M_n:=h^{\ssup1}(Y_n)\bigl(1-\lambda^{\ssup1}_S/4\bigr)^{-n}
\end{equation}
is a martingale with respect to the filtration $\sigma(Y_0,\dots,Y_n)$.

We will now derive the desired conclusion from this fact. We actually begin with an upper bound instead of a lower bound. Since $M_{n\wedge\tau_S}^2$ is a submartingale and $M_{n\wedge\tau_S}^2\le M_n^2$ due to the fact that $h^{\ssup1}(Y_{\tau_S})=0$ when $\tau_S<\infty$, we have
\begin{equation}
h^{\ssup1}(x)^2\le E^x\bigl(M_{n\wedge\tau_S}^2\bigr)\le E^x\bigl(M_n^2\bigr)
=\bigl(1-\lambda^{\ssup1}_S/4\bigr)^{-n}E^x\bigl(h^{\ssup1}(Y_n)^2\bigr).
\end{equation}
Now let~$S$ be as in the statement and set $n:=\lceil r(t,\beta)^2\rceil$. By \eqref{E:2.13ue}, the prefactor on the right is bounded uniformly in $t,\beta\ge1$. Since~$S$ fits into a ball of radius proportional to~$r(t,\beta)$, the Local Central Limit Theorem shows that $P^x(Y_n=y)\le cn^{-1}\le c r(t,\beta)^{-2}$ uniformly in~$x,y\in S$. Hence,
\begin{equation}
E^x\bigl(h^{\ssup1}(Y_n)^2\bigr)
=\sum_{y\in S}P^x(Y_n=y)h^{\ssup1}(y)^2\le cr(t,\beta)^{-2}\Vert h^{\ssup1}\Vert_{\ell^2(\BZ^2)}^2.
\end{equation}
As~$h^{\ssup1}$ is assumed normalized and $|S|$ is of order~$r(t,\beta)^2$, we get
\begin{equation}
\label{E:5.16a}
\max_{x\in S}\,h^{\ssup1}(x)^2\le\frac{\tilde c}{|S|}
\end{equation}
with~$\tilde c$ depending only on~$U$.

Moving over to the desired lower bound, the sequence $\{h^{\ssup1}(Y_{n\wedge\tau_S})\colon n\ge0\}$ is a non-negative supermartingale and so for any stopping time~$T$,
\begin{equation}
\label{E:5.17a}
h^{\ssup1}(0)\ge E^0\bigl(h^{\ssup1}(Y_{T\wedge\tau_S})\bigr).
\end{equation}
We will use this for~$T$ equal to $T_a:=\inf\{n\ge0\colon Y_n\in L_a\}$ where
\begin{equation}
L_a:=\Bigl\{x\in S\colon h^{\ssup1}(x)^2\ge\frac a{|S|}\Bigr\}.
\end{equation}
Then \eqref{E:5.17a} reads
\begin{equation}
h^{\ssup1}(0)\ge \frac{a}{|S|}P^0(T_a<\tau_S).
\end{equation}
It suffices to derive a lower bound on $P^0(T_a<\tau_S)$ for $S:=S((1-\epsilon)U,t,\beta)$ which is uniformly positive in the limit as~$r(t,\beta)\to\infty$ for all~$\epsilon$ small enough.

Since $h^{\ssup1}$ is normalized, the upper bound \eqref{E:5.16a} ensures that $|L_a|\ge\frac{1-a}{\tilde c-a}|S|$ (for $a<1$). We assume $a>0$ is small enough that $\frac{1-a}{\tilde c-a}>(2\tilde c)^{-1}$. Let $B_s(x):=x+[-s,s]^2\cap\BZ^2$. An averaging argument then shows that there is~$\delta>0$ depending only on~$\tilde c$ and~$U$ such that for at least one~$x\in S$,
\begin{equation}
B_{\delta r(t,\beta)}(x)\subset S\quad\text{and}\quad
\bigl|B_{\delta r(t,\beta)}(x)\cap L_a\bigr|\ge \frac1{2\tilde c}\bigl|B_{\delta r(t,\beta)}(x)\bigr|
\end{equation}
Now take $\tilde\delta:=\delta/(4\tilde c)$. A use of the pigeon-hole principle shows that, for some integer $n$ with $n/r(t,\beta)\in[2\tilde\delta,\delta]$ we then have
\begin{equation}
\bigl|\{z\in\BZ^2\colon |z-x|_\infty=n\}\cap L_a\bigr|\ge \frac n{4\tilde c}
\end{equation}
We now bound $P^0(T_a<\tau_S)$ by first requiring the walk to hit $B_{\tilde\delta r(t,\beta)}(x)$ before exiting from~$S$ and then to hit $\{z\in\BZ^2\colon |z-x|_\infty=n\}\cap L_a$ before exiting from $B_n(x)$. The former event has a uniformly positive probability due to the convergence of the walk to Brownian motion in the supremum norm and the fact that both $B_{\tilde\delta r(t,\beta)}(x)$ and~$S$ scale proportionally to~$r(t,\beta)$. The latter event has a uniformly positive probability by the fact that, for any~$N\ge1$ large enough and any $c>2$, the exit distributions from~$B_{cN}(x)$ for the walk started at~$x$ and the walk started at any~$z\in\partial B_N(x)$ are mutually absolutely continuous with bounds that are uniform in~$N$ and~$z$.

Having proved \eqref{E:3.19}, the proof of \eqref{E:3.17b} is now straightforward. Indeed, a successive use of Lemma~\ref{lemma-2.7} shows that
\begin{equation}
\frac{h^{\ssup1}(y)}{h^{\ssup1}(x)}\ge (2d)^{-\diam(S)},\qquad x,y\in S.
\end{equation}
Taking~$x:=0$, the lower bound in \eqref{E:3.19} now finishes the job.
\end{proofsect}

\section*{Acknowledgments}
\noindent
This research has been partially supported by NSF award DMS-1407558 and GA\v CR project P201/16-15238S.

\bibliographystyle{abbrv}

\end{document}